\documentstyle{article}
\input amssym.def
\input amssym.tex

\newcommand{\TG}{{\Tilde{G}}}

\newcommand{\hr}{{\widehat{r}}}
\newcommand{\hpi}{{\widehat{\pi}}}
\newcommand{\R}{{\Bbb R}}
\newcommand{\CP}{{\Bbb CP}}
\newcommand{\Tilde}{\widetilde}
\newcommand{\ov}{\overline}
\newcommand{\p}{{\partial}}
\newcommand{\PP}{{\bf P}}
\newcommand{\Symp}{{\rm Symp}}
\newcommand{\Diff}{{\rm Diff}}
\newcommand{\Map}{{\rm Map}}
\newcommand{\Proj}{{\rm Proj}}
\newcommand{\Fol}{{\rm Fol}}
\newcommand{\Fib}{{\rm Fib}}
\newcommand{\Dom}{{\rm Dom}}
\newcommand{\Aut}{{\rm Aut}}
\newcommand{\Ker}{{\rm Ker}}
\newcommand{\Ww}{{\cal W}}
\newcommand{\SO}{{\rm SO}}

\newcommand{\U}{{\rm U}}
\newcommand{\SU}{{\rm SU}}

\newcommand{\cp}{{{\CP}\,\!^2}}
\newcommand{\ocp}{{\ov{\CP}\,\!^2}}
\newcommand{\bcp}{{\ov{\CP}\,\!^2}}

\newcommand{\Aa}{{\cal A}}
\newcommand{\Mm}{{\cal M}}
\newcommand{\Ll}{{\cal L}}
\newcommand{\Jj}{{\cal J}}

\newcommand{\Dd}{{\cal D}}

\newcommand{\Nn}{{\cal N}}
\newcommand{\Ff}{{\cal F}}

\newcommand{\Vv}{{\cal V}}
\newcommand{\Ss}{{\cal S}}

\newcommand{\th}{{\theta}}
\newcommand{\Q}{{\Bbb Q}}
\newcommand{\N}{{\Bbb N}}
\newcommand{\Pp}{{\bf P}}
\newcommand{\T}{{\Bbb T}}
\newcommand{\Z}{{\Bbb Z}}
\newcommand{\C}{{\Bbb C}}

\newcommand{\al}{{\alpha}}

\newcommand{\be}{{\beta}}

\newcommand{\Om}{{\Omega}}
\newcommand{\om}{{\omega}}
\newcommand{\eps}{{\varepsilon}}
\newcommand{\de}{{\delta}}

\newcommand{\Ga}{{\Gamma}}
\newcommand{\ka}{{\kappa}}
\newcommand{\la}{{\lambda}}
\newcommand{\si}{{\sigma}}
\newcommand{\La}{{\Lambda}}

\newcommand{\io}{{\iota}}
\newcommand{\Si}{{\Sigma}}

\newcommand{\MS}{{\medskip}}

\newcommand{\NI}{{\noindent}}
\newcommand{\proof}[1]{\noindent{\bf Proof#1:\  }}

\newcommand{\QED}{\hfill {\bf QED} \medskip}

\newtheorem{theorem}{Theorem}[section]
\newtheorem{thm}[theorem]{Theorem}
\newtheorem{cor}[theorem]{Corollary}

\newtheorem{defn}[theorem]{Definition}

\newtheorem{remark}[theorem]{Remark}
\newtheorem{rmk}[theorem]{Remark}
\newtheorem{lemma}[theorem]{Lemma}

\newtheorem{prop}[theorem]{Proposition}
\newtheorem{proposition}[theorem]{Proposition}

\title{Topology of symplectomorphism groups of rational ruled surfaces}
\author{Miguel Abreu\thanks{Partially supported by NSF grant DMS 9304580,
while at the Institute for Advanced Study (1996/97), and  afterwards by 
FCT grant PCEX/C/MAT/44/96 and PRAXIS XXI through the Research Units Pluriannual 
Funding Program.} \\
Instituto Superior T\'ecnico, Lisbon, Portugal
\\ (mabreu@math.ist.utl.pt)
\and
Dusa McDuff\thanks{Partially
supported by NSF grant DMS 9704825.} \\ State University of New York
at Stony Brook, USA \\ (dusa@math.sunysb.edu)}

\date{October 6, 1999}

\begin{document}

\maketitle


\begin{abstract}
Let $M$ be either $S^2\times S^2$ or the one point blow-up $\cp\#\,\bcp$ 
of $\cp$.  In both cases $M$ carries a family of symplectic forms $\om_\la$,
where $\la > -1$ determines the cohomology class $[\om_\la]$.
This paper 
 calculates the rational (co)homology of the group $G_\la$ of
symplectomorphisms of $(M,\om_\la)$  as well as the rational
homotopy type of its classifying space $BG_\la$.  It turns out that each
group $G_\la$ contains a finite collection $K_k, k = 0,\dots,\ell = \ell(\la)$, 
of finite dimensional Lie
subgroups that generate its homotopy.  We show that 
these subgroups ``asymptotically commute", i.e. all the higher Whitehead
 products that they generate vanish as $\la\to \infty$.  However, for
each fixed $\la$ there is essentially one  nonvanishing product that gives
rise to a ``jumping generator" $w_\la$ in $H^*(G_\la)$ and to a single relation
in the rational cohomology ring $H^*(BG_\la)$.  An analog of this generator $w_\la$ 
was also seen by Kronheimer in his study of families of
symplectic forms on $4$-manifolds using Seiberg--Witten
theory.  Our methods involve a close study of the space of 
$\om_\la$-compatible almost complex structures on $M$. 
\end{abstract}

\section{Introduction}
\label{sec:intro}

\medskip

Rational ruled surfaces occur as
projectivizations $\PP(L_k\oplus\C)$, where $L_k$ is a complex line bundle
over $S^2$ with first Chern class $k$.
There are two cases to consider, firstly ($k$ is even)
when the
underlying manifold is diffeomorphic to the product $S^2\times S^2$ and
secondly ($k$ odd) when it is diffeomorphic to $\cp\#\, \ocp$, the one point
blow up of $\C P^2$.  Work of Taubes, Liu--Li and Lalonde--McDuff (for
detailed references see~\cite{LM}) implies that every symplectic 
form on one of these
manifolds is ``standard", i.e. it belongs to one of the two following families:
$$
M_\la^0 = (S^{2}\times S^{2}, \om^0_\la = 
(1+\la)\si_{0}\times\si_{0}),\quad  0\leq\la\in~
{\R},
$$
 where $\si_0$ is an area form on $S^2$ with total area $1$, and
$$
M_\la^1 =  (\cp\#\, \ocp, \om^1_\la),\quad -1 <\la\in\R,
$$
where $\om^1_\la$ takes the value $2 +\la $ on the class of the line
and $1 +\la $ on the exceptional divisor.  More explicitly, $\om_\la^1$ can
be obtained by collapsing the bounding spheres  of the annulus
$$
\{z \in \C^2: \la + 1\le \pi |z^2|\le \la + 2\}
$$
along their null foliations by Hopf circles.  We will denote
the corresponding groups of symplectomorphisms by $G_\la^i$, $i = 0,1$.

As was first observed by Gromov~\cite{G}, the
topological type of $G_\la^i$ changes as $\la$ increases.
He showed that the identity component of the group
$$
G_0^0 = \Symp(S^2\times S^2, \si_0\times \si_0)
$$
deformation retracts to the Lie group $\SO(3)\times \SO(3)$, and that
the full group is an extension of this Lie group by $\{\pm 1\}$, where $-1$
acts by permuting the factors.  He also pointed out that as soon as $\la$ gets
bigger than zero, $G_\la^0$ becomes connected and a new element of infinite
order appears in its fundamental group.  The key idea in his proof was to
look at the action of $G_\la^0$ on the contractible space $\Jj_\la^0$ of
$\om_\la^0$-compatible almost complex structures.  

By exploiting this idea further, the
 first author~\cite{Ab}  calculated the rational cohomology ring of
$G_\la^0$ for $\la$ in the range $0< \la \le 1$.  The first main result of the
present paper is the completion of this calculation for arbitrary $\la$.   Recall
that the  rational cohomology ring  $H^*(G)$ of a group $G$ is a free graded
ring. We denote by $\La(x_1,\dots, x_\ell)$ the exterior algebra over $\Q$ on
generators $x_j$ of odd degree and by $S(w_1,\dots,w_\ell)$ the  polynomial
algebra on even generators $w_j$.

\begin{thm}\label{thm:main0} When $\la > 0$ the group $G_\la^0$ is
path-connected, and  has fundamental group $\Z \oplus \Z/2\Z \oplus
\Z/2\Z$.  Moreover if
 $0\le \ell  -1 < \la\le \ell$ for some $\ell\in \N$,
$$
H^*(G_\la^0; \Q) = \La(a,x,y)\otimes S(w_\ell),
$$
where $\deg a = 1, \deg x = \deg y = 3$ and $\deg w_\ell = 4\ell$.
\end{thm}

Thus the cohomology remains stable as $\la$ varies within an interval of
the form $(\ell-1, \ell]$ but jumps as $\la$ moves past the endpoint $\ell$.
 The proof is simply a more
elaborate version of the calculation in~\cite{Ab}. It is based on the fact that
the space $\Jj_\la^0$  has a natural stratification by sets $U_j^0$, each of
which is homotopy equivalent to a homogeneous space of the group
$G_\la^0$.  More precisely, the stratum $U_k^0$  contains the Hirzebruch
integrable complex structure $J_{2k}$ with a holomorphic sphere of
self-intersection $-2k$, and it  is (weakly) homotopic to the quotient
$G_\la^0/\Aut(J_{2k})$, where $\Aut(J_{2k})$ denotes the stabilizer of
$J_{2k}$: see \S2. This gives rise to a family of Mayer--Vietoris sequences
which can be completely analyzed: see Proposition~\ref{prop:maps}.  
\MS

It is easiest to understand what is happening if we consider homotopy, rather
than cohomology.  We define $K_k^0$ to be 
the identity component of $\Aut(J_{2k})$.  Thus $K_0^0$ is the subgroup
$\SO(3)\times \SO(3)$ mentioned above, while, for $k > 0$, we will see in
\S2 that $K_k^0 $ is isomorphic to $S^1\times \SO(3)$.
In this setting,  the generators $\xi,\eta\in \pi_3(G_\la)$ dual to $x,y$ come from
the inclusion of $K_0^0$ into $G_\la^0$.  As $\la$ increases above $0$
a new $1$-dimensional generator $\al$ appears that lies in $\pi_1(K_1^0)$.
It does not commute with $\xi$, and so a new auxiliary class $\zeta_1$ also
appears that is represented in homotopy by the Samelson product
$[\al,\xi]$.\footnote{
This is given by the map $S^4 = S^1\times S^3/S^1\vee S^3
\to G_\la$ induced by the commutator 
$$
S^1\times S^3\to G_\la:\qquad (s,u)\mapsto \al(s) \xi(u) \al(s)^{-1} \xi(u)^{-1}.
$$
It is well known to be the desuspension of the Whitehead
product on $\pi_*(BG)$: see~\cite{Wh}.}
As $\la$ increases past $1$, one can find
representatives of $\al$ and $\xi$ that both lie in $K_2^0$.  Since they commute,
the Samelson product $[\al,\xi]$ vanishes.  Hence there is a  $5$-disk that
bounds $\zeta_1$, and the new $8$-dimensional generator $\zeta_2$ is a higher 
product made from this new disk and $\xi$.  It is shown in~\cite{Al2} that 
 such higher  products in the group $G$ are the desuspension of higher
Whitehead products in  $BG$ and hence give rise to relations in
the rational cohomology ring of the classifying space $BG$. 
By pursuing this argument, we prove the following result in \S\ref{ss:whit}.

\begin{theorem}\label{thm:bg}  If
 $0\le \ell  -1 < \la\le \ell$ for some $\ell\in \N$, then
$$
H^*(BG_\la^0; \Q) = S(A,X,Y)/\{A(X-Y)(4X - Y)\dots(\ell^2 X - Y) = 0\},
$$
where $\deg A = 2, \deg X = \deg Y = 4$.
\end{theorem}

A phenomenon similar to the existence of these ``jumping'' generators 
$w_{\ell}$ amd $\zeta_\ell$ was also
discovered by Kronheimer~\cite{K}  in the guise of some
nontrivial  families of symplectic forms, detected using 
properties of Seiberg--Witten invariants. The relation of our work to his is
explained in \S\ref{ss:kron}.
Kronheimer's paper was in turn motivated by work of Seidel~\cite{Sd} 
showing that
``many''  symplectic $4$-manifolds admit symplectic diffeomorphisms which
are  differentiably isotopic to the identity but not symplectically isotopic.
Of course, Theorems~\ref{thm:main0} and~\ref{thm:main1} 
show that this type of phenomenon does not 
happen on $S^2\times S^2$ or $\cp\#\, \ocp$. 

Another interesting question concerns the behavior of the groups $G_\la^0$
as $\la \to \infty$.   Since there is no obvious map $G_\la^0 \to G_\mu^0$
when $\la < \mu$ it is not quite clear how to interpret this limit.
 In \S\ref{sec:lim} we show that it can be defined as  a  bundle of
groups  $$
G_\infty^0 = \bigcup_{\la> 0} \{(\phi,\la): \phi\in G_\la^0\}\subset
\Diff(S^2\times S^2)\times \R^+. $$
One can think of $G_\infty^0$  as a topological category (or groupoid) with
$\R^+$ as space of objects and where the set of morphisms from $\la$ to
$\mu$ is empty unless $\la=\mu$, in which case it equals $G_\la^0$.
We  prove the following result.

\begin{prop}\label{prop:lim} $G_\infty^0$ is homotopy equivalent to the
group
$$
\SO(3)\times \Map(S^2,\SO(3)), $$
where $\Map$ denotes the space of  unbased smooth maps.
\end{prop}

Thus $G_\infty^0$ is homotopy equivalent to the group of all
diffeomorphisms of $S^2\times S^2$ that commute with the projection to the
first (i.e. the larger)  factor.  In other words, in the limit the only remaining structure
is that of the fibration $S^2\times S^2\to S^2$.  Note that in this limit we allowed the
size of the larger sphere to increase without bound, while the smaller sphere has
area fixed at $1$.  Sometimes it is convenient to reparametrize this, fixing
the area of the larger sphere at $1$ and allowing the size of  the smaller
sphere to go to zero.  (This is the kind of adiabatic limit considered by
Salamon in~\cite{Sa} for example.)  Since the group of symplectomorphisms 
does  not change when the symplectic form is multiplied by a constant,
these different ways of taking the limit yield the same result.  
Observe also that an analogous result holds for the symplectomorphism
groups of any ruled surface: see Remark~\ref{rmk:any}. \MS

There is a completely analogous story for the groups $G_\la^1$ of
symplectomorphisms of $\cp\#\,\ocp$.  In this case Gromov's methods
 show that $G_\la^1$ deformation retracts onto the
unitary group $\U(2)$ when $-1 < \la \le 0$, and our arguments again show that
its cohomology jumps as $\la$ passes each integer.  This time, however,
there already are generators  $a ,x$ of degrees $1$ and $3$ in the
cohomology of $\U(2)$ and  the new persistent generator $y$ that appears
has degree $3$, while the ``jumping generator" has degree $6, 10,\dots$ and
so on.  

\begin{thm}\label{thm:main1} For all $\la > -1$ the group $G_\la^1$ is
path-connected, and  has fundamental group $\Z$.  Moreover if
 $0\le \ell-1  < \la\le \ell$ for some $\ell\in \N$,
$$
H^*(G_\la^1; \Q) = \La(a,x,y)\otimes S(w_\ell),
$$
where $\deg a = 1, \deg x = \deg y = 3$ and $\deg w_\ell = 4\ell + 2$.
\end{thm}

\begin{theorem}\label{thm:bg1}  If
 $0\le \ell  -1 < \la\le \ell$ for some $\ell\in \N$, then
$$
H^*(BG_\la^1; \Q) = S(A,X,Y)/\{(X-Y)(4X-Y)\dots((\ell+1)^2 X - Y) = 0\},
$$
where $\deg A = 2, \deg X = \deg Y = 4$.
\end{theorem}

\MS

\MS

\noindent {\bf Acknowledgements:} Part of this
work was carried out while the first author was a
postdoctoral member of the School of Mathematics of the Institute for
Advanced Study. He takes this opportunity to thank the staff and faculty of
the Institute for providing a wonderful environment, both for him and his
family.  The second author would like to thank Allday for bringing the 
references~\cite{Al,Al2,AA} to her attention.

\tableofcontents

\section{The stratification of $\Jj_\la^i$ and the groups $K_k^i$}

In this section we define the stratification of the space  $\Jj_\la^i$ of
$\om_\la^i$-compatible almost complex structures 
and establish its basic properties.  Each stratum is homotopy equivalent to a
homogeneous space $G_\la^i/K_k^i$, where the subgroups $K_k^i$
are finite dimensional Lie groups.  In \S\ref{ss:Kk} we work out the
relations between these groups $K_k^i$.

\subsection{Structure of $J$-holomorphic spheres on $M^i_\la$}
\label{structure}

We first consider  $M^0_{\la}= (S^{2}\times
S^{2}, \omega^0_{\la})$.  In this case we always take $\ell$ to be the
positive integer such that
$$
\ell -1<\la\leq \ell.
$$
Further,  $E_0$ denotes the homology class of $S^{2}\times\{pt\}$, $F$
denotes the fiber class $\{pt\}\times S^{2}$, and we define
$E_{2k}=E_0 - kF$ for $k=0,1,2,\ldots$. General background information about
$J$-holomorphic spheres can be found in~\cite{Ab,LM,MS1,MS2}.

Let $J_0=j_0\oplus j_0 \in \Jj^0_\la$ be the standard split compatible
complex structure. For $J_0$, the classes $E_0$ and $F$ are both
represented by $2$-parameter families of  holomorphic spheres
given respectively by $S^{2}\times\{s\}$ and $\{s\}\times S^{2}$, for any
$s\in S^{2}$.   Moreover,  by positivity of intersections no class $E_{2k}$ with
$k> 0$ can be represented by a $J_0$-holomorphic curve.  More
generally, we have:

\begin{proposition} \label{prop:struct0}
\begin{itemize}
\item[{\rm (i)}] 
 For any $J\in\Jj^0_\la$, the fiber class $F$ is represented by a
$2$-parameter family of embedded $J$-holomorphic spheres that
form the fibers of a fibration $S^2\times S^2\to S^2$. If $\la=0$ the same is
true for $E_0$. 

\item[{\rm (ii)}] For each $k, 0\le k\le \ell,$ the class
$E_{2k}$ is represented by a $J$-holomorphic sphere
for some $J\in \Jj_\la^0$. 

\item[{\rm (iii)}] Conversely, for each $J\in \Jj_\la^0$,
there is an integer $0\le  k\leq \ell$ such that  the class $E_{2k}$
has a  $J$-holomorphic representative.   When $k > 0$ this sphere is unique
and no other class $E_{2j}$ with $0\le j, j\neq k$, has such a
representative.
\end{itemize}
 \end{proposition}

\proof{} The first  statement is a well-known and standard result 
in the theory of  $J$-holomorphic curves: see e.g.~\cite{LM}.
To prove (ii), first note that if $0\le k\leq \ell$, then
$$
\om^0_\la(E_{2k})=1+\la -k > 0,
$$
which is a necessary condition for $E_{2k} = E_0 - kF$ to have a holomorphic
representative for some $J\in \Jj_\la^0$.  We show that it does  in
\S\ref{ss:red} by an explicit construction.

The third statement follows from the compactness theorem for
$J$-holomorphic curves.  Since $\Jj_\la^0$ is connected there is a path $J_t,
t\in [0,1]$ connecting the standard split structure $J_0$ to any given $J = J_1$.  For a
generic path, the class $E_0$ will be represented for all $J_t, t < 1,$ and will
degenerate at $t=1$ into some cusp-curve that has to contain a component in some
class $E_{2k}$.  
 Next recall that the principle of 
positivity of intersections says that each
intersection point of two distinct $J$-holomorphic curves in a
$4$-manifold contributes positively to their intersection number.  This
implies that the intersection number $A\cdot B$ of any
 two different classes $A,B$ with $J$-holomorphic
representatives  must be $\ge 0$.  Since
$$ 
E_{2k}\cdot E_{2j} = (E_0 - kF)\cdot(E_0 - jF) = -k-j,
$$
the uniqueness and nonexistence statements in (iii)  follow
immediately. \QED

Consider now $M^1_\la = (\cp\#\,\ocp, \om^1_\la)$, with $-1 <\la\in\R$,
described in \S\ref{sec:intro}. Recall that $E_1$ denotes the exceptional
divisor,
with self-intersection $-1$, and $F$ is used again to denote the fiber
of the fibration $\cp\#\,\ocp\rightarrow S^2$, having self-intersection
$0$. Moreover, the symplectic form $\om^1_\la$ is such that
$$
\om^1_\la(F)= 1\ {\rm and}\ \om^1_\la(E_1)=\la + 1\ .$$
  In this case also we define the integer $\ell$ by the
condition $$
\ell -1 < \la\le \ell.
$$
We further  define $E_{2k+1}$ by $E_{2k+1}=E_1 - kF$ for $k=0,1,2,\ldots\ $.
Thus, again, the self-intersection number of $E_{2k+1}$ is $2k+1$.

Let $J_1\in\Jj^1_\la$ be the standard compatible complex structure on
$\cp\#\,\ocp$. For $J_1$, the class $F$ is represented by a $2$-parameter
family of holomorphic spheres given by the fibers of the holomorphic
fibration $\cp\#\,\ocp\rightarrow S^2$, and the class $E_1$ is represented
by a unique holomorphic sphere given by the
exceptional divisor.   Here is the analog of Proposition~\ref{prop:struct0}. 
The proof is similar and is left to the reader.

\begin{proposition} \label{prop:struct1}
\begin{itemize}
\item[{\rm (i)}] 
  For any $J\in\Jj^1_\la$, the fiber $F$ is represented by a
$2$-parameter family of embedded $J$-holomorphic spheres that
fiber $\cp\#\,\ocp$. If $\ell = 0$, $E_1$  is
represented by a unique embedded $J$-holomorphic sphere
 for any  $J\in\Jj^1_\la$.

\item[{\rm (ii)}]   For each $k, 0\le k\le \ell,$ the class
$E_{2k+1}$ is represented by a $J$-holomorphic sphere
for some $J\in \Jj_\la^1$. 

\item[{\rm (iii)}]  Conversely, for each $J\in \Jj_\la^1$,
there is an integer $ k\leq \ell$ such that  the class $E_{2k+1}$
has a  $J$-holomorphic representative.   This sphere is unique
and no other class $E_{2j+1}$ with $0\le j, j\neq k$, has such a
representative.
\end{itemize}
\end{proposition}

\subsection{The strata of $\Jj^i_\la$} \label{strata}

As always we fix $\ell\in\N$ and assume that
$\ell-1<\la\leq \ell$. For each $0\leq k\leq \ell$ and $i = 0,1$, define
$U^i_k\subset\Jj^i_\la$ as the set of all $\om^i_\la$-compatible
almost complex structures for which the class $E_i - kF = E_{2k+i}$ 
is represented by a pseudo-holomorphic sphere.   A consequence of
Propositions~\ref{prop:struct0} and~\ref{prop:struct1} is that the
$U^i_k$'s  are disjoint subsets of $\Jj^i_\la$ for $i = 0,1$.
  Throughout we  work with
 $C^\infty$-maps and almost complex structures, and so by manifold we mean a
Fr\'echet manifold.

\begin{proposition} \label{prop:strata} Let $i = 0$ or $1$.  Then
$U^i_0$ is an open dense subset of $\Jj^i_\la$, and
each $U^i_k$, $1\leq k\leq
\ell$, is a non-empty submanifold of
$\Jj^i_\la$ of codimension $4k-2$ if $i = 0$ and  $4k$ if $i = 1$. Its
 closure is given by 
$$\overline{U}^i_k = U^i_k \cup U^i_{k+1}
\cup\cdots\cup U^i_\ell.
$$
In particular,
$$
\Jj^i_\la = \overline{U}^i_0 = U^i_0 \cup U^i_{1} \cup\cdots\cup U^i_\ell. $$
\end{proposition}
\proof{} This is a relatively easy result whose proof is analogous to the proof
of Theorem 1.8 in~\cite{Ab}.  See also \S4.1 in~\cite{M5}. \QED

For the algebraic computations in \S\ref{sec:algebra}, we
need more information regarding the way the $U^i_k$'s 
fit together to give the contractible space $\Jj^i_\la$. Define $U^i_{01\ldots
k}$ for $0\le k \le \ell$  by
$$
U^i_{01\ldots k} = U^i_0 \cup U^i_{1} \cup\cdots\cup U^i_k.
$$
The next result says that these spaces $U$ are stratified.  By this 
we mean  that  $U$ is a union of a finite number of disjoint manifolds that are called
strata.  Each stratum $\Ss$ has a neighborhood $\Nn_\Ss$ that projects to $\Ss$ by a
map $\Nn_\Ss\to \Ss$.  When $\Nn_\Ss$ is given the
induced stratification, this map is a locally trivial
fiber bundle whose fiber has the
form of a cone $C(\Ll)$  over a finite dimensional stratified space $\Ll$ that is
called the {\em link} of $\Ss$ in  $U$.   Moreover, $\Ss$ sits inside $\Nn_\Ss$
as the set of vertices of all these cones.

\begin{proposition} \label{prop:link}
\begin{itemize}
\item[{\rm (i)}]  For each $0\leq k\leq \ell$, the space $U^i_{01\ldots k}$  is stratified with
strata $U^i_j$, $0\leq j\leq k$. 

\item[{\rm (ii)}]  The link of $U^i_k$  in $U^i_{01\ldots k}$
 is the sphere $S^{4k-3}$ when $i=0$ and $S^{4k-1}$ when $i = 1$. 
Moreover,  if $k>1$ the link $L^i_k$  of $U^i_k$ in
$U^i_{k-1}$  is a smooth $3$-manifold (in fact
a lens space).
\end{itemize}
 \end{proposition}
\proof{} The first statement is the main result of~\cite{M5} and is proved by 
showing that  for all $J\in \Jj_\la^i$ the space of all $J$-holomorphic stable
maps into $M_\la^i$ has a corresponding stratification.  
The first statement in (ii) follows from the previous proposition, while
the second  can
be proved by looking at the gluing parameters necessary to glue in one fiber:
see~\cite{M5}.  We  give a rough idea of its proof since this will be
needed later.

  Each element $J\in U_k$ defines a fibration
$\pi_J: M_\la^i \to S^2$ together with a unique $J$-holomorphic section
$C_J$ in the class $E_{2k+i} = E_i - kF$.  There is a corresponding family 
$\tau_q, q\in S^2$, of $J$-holomorphic cusp-curves (or stable maps) in class
$E_i-(k-1)F$.  Namely, $\tau_q$ is the union of $C_J$ with the fiber
$F_q= \pi_J^{-1}(q)$.  The gluing process allows one to build a family of
$J_{(q,a)}$-holomorphic spheres $S_{(q,a)}$ from the pair of components in $\tau_q$,
one for each sufficiently small ``gluing parameter" $a\in \C_q$. 
Here $\C_q$ can be identified with the tensor product of the tangent spaces
to the two components $C_J$ and $F_x$ of $\tau_q$ at their
intersection point, and so fit together to
form a complex line bundle $L$  over $S^2$ with Chern class $2 - 2k$.  The
permitted gluing parameters $(q,a)$ fill out a neighborhood $\Vv$ of its zero
section.

Now observe that because
$S_{(q,a)}$ lies in the class $E_i-(k-1)F$, the corresponding almost complex
structures $J_{(q,a)}$ lie in $U_{k-1}$.  Further, except in the case $i = 0,
k = 1$, there is a unique $J_{(q,a)}$-holomorphic curve in the class
$E_i-(k-1)F$, and one can show that the map
$$
\Vv\to U_{k-1}:\quad {(q,a)}\mapsto J_{(q,a)}
$$  
is a homeomorphism onto  a normal slice to $U_k$ in
$U_{k-1}$ at $J$.  Thus the link  $L_k^i$ is homeomorphic to
$\p\Vv$ which is the unit sphere bundle of $L$ and so is a  lens space.
 \QED

\subsection{The manifolds $M^i_\la$ as symplectic reductions
of $\C^4$} \label{ss:red}

 In order to be able to identify the strata $U^i_k$ of
$\Jj^i_\la$ with quotients of $G^i_\la$ by appropriate compact subgroups, we now
show how  to describe $M^i_\la$ as a symplectic
reduction of $\C^4$. 
A good reference for this subsection is the article by M.Audin in~\cite{AL}.
For general expositions on symplectic group actions and reduction see,
for example,~\cite{A} and~\cite{MS1}.

Consider the unitary action of ${\T}^{2}$
on ${\C}^{4}$ given by
$$
(s,t)\cdot(z_{1},z_{2},z_{3},z_{4})=(s^{m}tz_{1},tz_{2},sz_{3},sz_{4}),
\ |s|=|t|=1,
$$
and let $\phi :{\C}^{4}\rightarrow{\R}^{2}$ be the corresponding
moment map:
$$
\phi(z_{1},z_{2},z_{3},z_{4})=\frac{1}{2}(m|z_{1}|^{2}+|z_{3}|^{2}+|z_{4}|^{2},
|z_{1}|^{2}+|z_{2}|^{2}).
$$
If $\mu>m$ one checks easily that $\xi=(\mu,1)\in{\R}^{2}$ is a
regular value of $\phi$ and so we can consider the symplectic reduction
$\phi^{-1}(\mu,1)/{\T}^{2}$ which we denote by $R^m_{\mu}$.

Since the reduction is performed with respect to a unitary action of
$\T^2$, the complex structure of $\C^4$ descends naturally to the quotient,
giving $R^m_\mu$ the structure of a K\"ahler manifold.
By looking at the action of $T^2$ on the level set
$\phi^{-1}(\mu, 1)$, it is not hard to see that
 $R^{m}_\mu$ can be identified with the space
$$
W_m= \{([u_0, u_1], [w_0: w_1: w_2]) : u_0^m w_1 = u_1^m w_0\}\subset
\CP^1\times \CP^2
$$
via the map 
$$
\C^4\to W_m: (z_1,z_2,z_3,z_4)\mapsto ([z_3: z_4], [z_3^m z_2: z_4^m z_2:
z_1]).
$$
(See also~\cite{AL} page 61.) Hence $W_m$ fibers over the copy of $S^2$
with coordinates $z_3, z_4$, with fiber the $z_1,z_2$-sphere.
Standard arguments show that when $m$ is even this fiber bundle is trivial, 
and that when $m$ is odd it is nontrivial.

Next observe that the
submanifolds $\Tilde{E}_m\subset R^m_\mu$ defined by
$$
\Tilde{E}_m=\{ {\rm\ reduction\ of\ } \{z_{2}=0\}\subset{\C}^{4}\}$$
are holomorphic spheres with symplectic area given by $\mu - m$.
Thus, when $m=2k$ is even,  $R^{2k}_\mu$ 
is diffeomorphic to 
$S^2\times S^2$ 
 and the submanifold $\Tilde{E}_{2k}$ represents the class 
$E_{2k} = E_0 -kF$ defined in \S\ref{structure}. Moreover, if 
$\mu = 1+\la +k >  2k$ we have that $R^{2k}_{1+\la +k}$
is symplectomorphic to $M^0_\la$ and, if we define $J_{2k}$ as the push-forward
of the complex structure of $R^{2k}_{1+\la +k}$ under this symplectomorphism, 
the class $E_{2k}$ has a $J_{2k}$-holomorphic representative given by the 
image of $\Tilde{E}_{2k}$, i.e. $J_{2k}\in U^0_k$. Note that $J_0$ defined this 
way is the standard split complex structure $j_0\oplus j_0$ considered in 
\S\ref{structure}.

Similarly, when $m=2k+1$ is odd, $R^{2k+1}_\mu$ is
diffeomorphic to $\cp\#\ocp$ and the submanifold $\Tilde{E}_{2k+1}$
represents the class $E_{2k+1} = E_1 -kF$ also defined
in \S\ref{structure}. If $\mu = 2+\la +k > 2k+1$ we also have that
$R^{2k+1}_{2+\la +k}$
is symplectomorphic to $M^1_\la$ and, if we define $J_{2k+1}$ as the 
push-forward of the complex structure of $R^{2k+1}_{2+\la +k}$ under 
this symplectomorphism, the class $E_{2k+1}$ has a $J_{2k+1}$-holomorphic 
representative given by the image of 
$\Tilde{E}_{2k+1}$, i.e. $J_{2k+1}\in U^1_k$.\MS
\MS

We now want to understand the isometries of the above K\"ahler structures on
$M^0_\la$ and $M^1_\la$. When $m=0$, the subgroup $\SU(2)\times
\SU(2)\subset \U(4)$ \ given by matrices of the form
\[ \left( \begin{array}{cc} A & 0 \\ 0 & D \end{array} \right) {\rm \ with\ }
A,D\in SU(2), \]
acts on ${\C}^{4}$ by matrix multiplication, preserving the moment map $\phi$ and
commuting with the $\T^2$ action defined above.
Since the intersection of these two subgroups of $\U(4)$ is
$\{\pm 1\}\times\{\pm 1\}$, we get an effective K\"ahler action of
$\SO(3)\times \SO(3)\cong \SU(2)/\{\pm 1\}\times \SU(2)/\{\pm 1\}$ on
$R^0_{\mu}$, which is nothing else than the standard action of $\SO(3)\times
\SO(3)$ on $S^2\times S^2$ (see~\cite{I}).

When $m>0$, the above construction can still be done if we restrict ourselves to
the subgroup $\SU(2)\times S^1 \subset \U(4)$ given by matrices of the form
\[ 
\left( \begin{array}{cc} A & 0 \\ 0 & D \end{array} \right) {\rm \ with\ }
A=\left( \begin{array}{cc} \alpha & 0 \\ 0 & \overline{\alpha} \end{array}
\right) , \qquad\alpha\overline{\alpha}=1, {\rm \ and\ } D\in \SU(2). 
\]
This subgroup also acts on ${\C}^{4}$ by matrix multiplication, preserving
the moment map $\phi$ and commuting with the above $\T^2$ action.

When $m=2k$ is even, the intersection of these two subgroups of $\U(4)$ is
again $\{\pm 1\}\times\{\pm 1\}$ and we get an effective K\"ahler action of
$\SO(3)\times S^1 \cong \SU(2)/\{\pm 1\}\times S^1 /\{\pm 1\}$ on
$R^{2k}_{\mu}$. Moreover, the submanifold $\tilde{E}_{2k}\subset
R^{2k}_{\mu}$ is a connected component of the fixed point set of the $S^1$
part of this action. Taking $\mu = 1 + \la + k$ and using the fact that
$R^{2k}_{1+\la +k}$ is symplectomorphic to $M^{0}_{\la}$, we have
constructed this way a $\SO(3)\times S^1$ K\"ahler action on $M^{0}_{\la}$,
with complex structure $J_{2k}\in U^{0}_{k}$, such that each point of the
$J_{2k}$-holomorphic representative for the class $E_{2k}$ is fixed by the
$S^1$ part of the action. Note that when $k=1$ the $\SO(3)$ part of this action
is the same as the diagonal $\SO(3)$ action on $S^2\times S^2$
(see~\cite{Ab}).

When $m=2k+1$ is odd, the intersection of the above $\SU(2)\times S^1
\subset \U(4)$ and $\T^2\subset \U(4)$ is given by $\{ (1,\pm 1),(-1,\pm i)
\}$, and so we get an effective K\"ahler action of
$\U(2)\cong \SU(2)\times S^1 / \{ (1,\pm 1),(-1,\pm i) \}$ on
$R^{2k+1}_{\mu}$. Moreover, the submanifold $\tilde{E}_{2k+1}\subset
R^{2k+1}_{\mu}$ is a connected component of the fixed point set of the
natural $S^1$-subgroup that this description of $\U(2)$ gives:
$$ 
S^1\cong S^1/\{\pm 1\} \hookrightarrow \{ 1\}\times S^1 / \{\pm 1\}
\hookrightarrow \SU(2)\times S^1 / \{ (1,\pm 1),(-1,\pm i) \}\ .
$$
Taking $\mu = 2 + \la + k$ and using the fact that $R^{2k+1}_{2+\la +k}$
is symplectomorphic to $M^{1}_{\la}$, we have constructed this way a
$\U(2)$ K\"ahler action on $M^1_\la$, with complex structure $J_{2k+1}\in
U^{1}_{k}$, such that the $J_{2k+1}$-holomorphic representative for the
class $E_{2k+1}$ is fixed by a given $S^1$ subgroup of $\U(2)$.

We have thus proved the following:

\begin{prop} \label{prop:reduction}
\begin{itemize}
\item[(i)] For each $0< k\leq \ell$, $M^0_\la$ has a symplectic
$K^0_k \cong \SO(3)\times S^1$ action, which is K\"ahler for a standard
complex structure $J_{2k}\in U^0_k$ and is such that the
$J_{2k}$-holomorphic representative $C_{2k}$ for the class $E_{2k}$ is fixed
by the $S^1$ part of the action.
\item[(ii)] For each $0\leq k\leq \ell$, $M^1_\la$ has a symplectic
$K^1_k \cong \U(2)$ action, which is K\"ahler for a standard
complex structure $J_{2k+1}\in U^1_k$ and such that the $J_{2k+1}$-holomorphic
representative $C_{2k+1}$ for the class $E_{2k+1}$ is fixed by a given $S^1$ 
subgroup of $\U(2)$.
\end{itemize}
\end{prop}

\subsection{Geometric description of the strata}
\label{description}

The next task is to identify the stratum $U^i_k$ of $\Jj^i_\la$ with the
quotient of $G^i_\la$ by the isometry group $K_k^i$.

\begin{prop} \label{prop:description}
\begin{itemize}
\item[{\rm (i)}] If $\la > 0$, the stratum $U^0_0\subset\Jj^0_\la$ is weakly homotopy
equivalent to $G^0_\la / SO(3)\times SO(3)$.

\item[{\rm (ii)}] Each strata $U^0_k$ of $\Jj^0_\la$, $1\leq k\leq \ell$, is weakly
homotopy equivalent to $G^0_\la / K^0_k$, where $K^0_k \cong SO(3)\times S^1$
is the group of K\"ahler isometries of a standard complex structure
$J_{2k}\in U^0_k$.

\item[{\rm (iii)}] Each strata $U^1_k$ of $\Jj^1_\la$, $0\leq k\leq \ell$, is weakly
homotopy equivalent to $G^1_\la / K^1_k$, where $K^1_k \cong U(2)$
is the group of K\"ahler isometries of a standard complex structure
$J_{2k+1}\in U^1_k$.
\end{itemize}
\end{prop}
\proof{} Part (i)  was  proved in~\cite{Ab}. When $\la = 0$ the corresponding
statement was established by Gromov in~\cite{G}.  Here the group $K_0^0 =
\SO(3)\times \SO(3)$ is not the full isometry group of $J_0$: one has to add in the
involution that interchanges the two sphere factors.  Gromov then showed that 
the quotient of $G_0^0$ by this isometry group is $U_0^0$, which is this case is
contractible since it is the whole of $\Jj_0^0$.

The case $k=1$ of part (ii) was proved in~\cite{Ab}. Since that
proof generalizes directly to any $k$, with $1\leq k\leq \ell$, we
just recall its main steps adapted to this more general context.

One denotes by $\Ss^0_k$ the space of embedded symplectic $2$-spheres
in $M^0_{\la}$ representing the homology class $E_{2k}$. Associated
to any $J\in U^0_k$ we have a well defined element of
$\Ss^0_k$ given by the unique $J$-sphere representing
$E_{2k}$. Conversely, given any element in $\Ss^0_k$
the space of all $J\in U^0_k$ that make it $J$-holomorphic is
non-empty and contractible. It follows that $U^0_k$ is
weakly homotopy equivalent to $\Ss^0_k$.

The next step is to prove that $G^0_\la$ acts transitively on
$\Ss^0_k$.   Here is an outline of the method.  First of all, given two curves $C_j,
j = 1,2,$ in $\Ss^0_k$,  choose elements $J_j$ in $\Jj_\la^0$ so that $C_j$ is
$J_j$-holomorphic.  Then $C_j$ is 
a section
of the fibration of $M_\la^0$ formed by the $J_j$-holomorphic spheres in class $F$. 
It follows that there is a diffeomorphism $\phi$ of $M_\la^0$ that takes
$C_1$ to $C_2$.  It is easy to adjust $\phi$ so that it is a symplectomorphism near
$C_1$, and it can be adjusted to be a symplectomorphism everywhere because 
any two symplectic forms on the complement $M_\la^0 - C_1$ that are standard
near infinity are symplectomorphic.   The last step can be justified by 
explicit geometric
arguments as in~\cite{Ab} or by following the ideas
in~\cite{M2}.

Denoting by $H^0_k$ the subgroup of $G^0_\la$ consisting of
symplectomorphisms that preserve the unique $J_{2k}$-holomorphic
sphere  in the class $E_{2k}$, we then have that
$G^0_\la / H^0_\la \cong \Ss^0_k$. The obvious map
$$ G^0_\la / K^0_k \rightarrow G^0_\la / H^0_k $$
is a fibration, with fiber over the identity given by
$H^0_k / K^0_k$. The last step of the proof of (ii) is to show
that $H^0_k / K^0_k$ is weakly contractible.
As in the proof of Lemma~3.6 in~\cite{Ab} for the case $k=1$,
one proves first that $H^0_k / K^0_k$ is homotopy equivalent to the 
subgroup of $G^0_\la$ consisting of symplectomorphisms that restrict 
to the identity on a neighborhood of the sphere $C_{2k}$, and then
shows that this last subgroup is weakly contractible.

Part (iii) of the theorem follows from the same steps as part (ii).
\QED

We are now in a position to calculate $\pi_0$ and $\pi_1$ of the groups 
$G_\la^i$. Since we already know the homotopy types of $G_0^0$ and 
$G_\la^1,\ -1 < \la \le 0$ (namely, an extension of $\SO(3)\times\SO(3)$ and  
$U(2)$), we will concentrate on the other cases.

\begin{cor}\label{cor:desc}
\begin{itemize}
\item[{\rm (i)}] The group $G_\la^i$ is connected, except when $i = \la = 0$.

\item[{\rm (ii)}] When $\la > 0$, $\pi_1(G_\la^0) = \Z\oplus\Z/2\Z\oplus \Z/2\Z$.

\item[{\rm (iii)}] When $\la > 0$, $\pi_1(G_\la^1) = \Z$.
\end{itemize}\end{cor}
\proof{}  The previous proposition implies that for $\la > 0$  the sequence
$$
K_k^i \to G_\la^i \to U_k^i
$$
is a homotopy fibration, in the sense that there is a corresponding long exact
sequence of homotopy groups.  Further $\pi_0(K_0^i) = 0$.

Since the strata $U_k^i$, $k > 0$, have codimension at least $2$ in $\Jj_\la^i$,
the top stratum $U_0^i$ is always path connected.  Hence $\pi_0(G_\la^i) =
\pi_0(K_0^i) = 0$.  Similarly, when $i=1$, the fact that the strata $U_k^1, k>0,$ have
codimension at least $4$ in $\Jj_\la^1$ implies that $U_0^1$ is $2$-connected.
Hence 
$$
\pi_1(G_\la^1) = \pi_1(K_0^1)=  \pi_1(U(2)) = \Z.
$$
However, the stratum $U_1^0$ has codimension $2$ in $\Jj_\la^0$ which means that 
there is a potential generator of $\pi_1(U_0^0)$ given by the loop which circles 
once around $U_1^0$.  The argument in~\cite{Ab} which shows that this loop
generates $\pi_1(U_0^0)$ when $\la \le 1$  goes through in the case $\la > 1$ 
without change.  
The rest of (ii) also follows  as
in~\cite{Ab}.\QED

\begin{cor}  Each  stratum $U_k^i$ is connected.
\end{cor}
\proof{}  This follows from Proposition~\ref{prop:description} and the fact
that  the groups $G_\la^i$ are connected except when $i = \la = 0$.
\QED

The group $K_k^i $ is defined to be the identity component of the stabilizer of
the element $J_{2k+i}$ of the stratum $U_k^i$.   We now show that it acts
transitively on the link  $L_k^i$  of $U_{k}^i$ in $U_{k-1}^i$ at the point
$J_{2k+i}$.  The case $i=0, k=1$ is special, since here the link is a circle.  In
all other cases the link has dimension $3$.

\begin{prop}\label{prop:link2}  When $k + i> 1$ the group
$K_k^i$ acts transitively on the link  $L_k^i$ of  $U_{k}^i$ in $U_{k-1}^i$
at $J = J_{2k+i}$.
\end{prop}
\proof{}  This is immediate from the description of this link given in
Proposition~\ref{prop:link}.   The group $K_k^i$ acts
transitively on the  sphere $C_J$ (where $J ={J_{2k+i}}$) and on
its unit tangent bundle.  It also acts transitively on the unit
tangent bundle to the fibers of $\pi_J$ 
at their point of intersection with $C_{J}$.  
The bundle
$L$ of gluing parameters is the tensor product of these two bundles, and it
is not hard to see that the action there is transitive too.
\QED

\subsection{Relations between the $K_k^i$}\label{ss:Kk}

We now discuss the relations between the elements in
$\pi_1(G_\la^i)$ and $H_*(G_\la^i)$ represented by the homotopy
and homology of the different subgroups $K_k^i$ in  $G_\la^i$, beginning as
usual with the case $i=0$.  Our arguments use some results from~\S\ref{ss:0}
and~\S\ref{sec:lim}.  However, the conclusions will not 
be used until we calculate  $H^*(BG_\la)$ in~\S\ref{ss:whit}.

It will be convenient to use the fact, proved in
\S\ref{sec:lim}, that the $G_\la^0$  have a limit at infinity $G_\infty^0$ that is
homotopy equivalent to the group $\Dd = \Dd^0$  of orientation preserving 
fiberwise diffeomorphisms of $S^2\times S^2$.  Since $\Diff^+(S^2)$
deformation retracts to $\SO(3)$, $\Dd$ is homotopy equivalent to the
product 
$$
\Dd = \SO(3)\times \Map(S^2, \SO(3)),
$$
where $\Map(Y,Z)$ denotes the space of smooth maps from $Y$ to $Z$.
 It is easy to see from the explicit description given in \S\ref{ss:red}
of the action of $K_k$ on $W_{2k}$ that this action commutes 
with the projection onto
the $z_3, z_4$ coordinates, i.e  with the projection of $S^2\times S^2$ onto its
first  factor.
Thus each $K_k^0$ is a subgroup of $\Dd$.

 The group $K_0^0 $ is simply $
\SO(3)\times \SO(3)$ and we will denote   the two obvious generators of
$\pi_1(K_0^0)$ by $\tau$ and $\tau'$.  For $k > 0$   let $\al_k$ be the
element of $\pi_1(K_k^0)$ corresponding to the $S^1$ action that fixes  the
sphere $C_{2k}$.  Further, let $\xi = \xi_0$ and $\eta$ be the two obvious 
spherical generators of $H_3(K_0^0, \Z)$ and  $\xi_k$,  for $k > 0$, be the
spherical generator of  $H_3(K_k^0, \Z)$.  Thus these classes are represented
by suitable images of $\SU(2) = S^3$ in the groups $K_k^0$.

\begin{lemma}\label{le:htpyK0}
 For each $k > 0$,
 \begin{itemize}
\item[{\rm (i)}]
$
\al_{k} = k \al_1 + (k-1) \tau' \in \pi_1(G_\la),
$
\item[{\rm (ii)}] $  \xi_{k} = \xi_0 + k^2\eta \in H_3(G_\la; \Q).$
\end{itemize}
\end{lemma}
\proof{} (i)   We saw in Corollary~\ref{cor:desc} that 
$\pi_1(G_\la^0) = \Z \oplus \Z/2\Z \oplus \Z/2\Z$, with generators lying in
$\pi_1(K_{k}^0)$ for $k = 0,1$.   We claim that
 $\pi_1(\Dd)$ is also
equal to $ \Z \oplus \Z/2\Z \oplus \Z/2\Z$ and has the same generators.
To see this, note first that $\tau$ generates the copy of $\Z/2\Z$ coming
from the base, i.e. 
from the fundamental group of the $\SO(3)$-factor in $\Dd$.
Elements of $\pi_1(\Map(S^2,\SO(3)))$ correspond to homotopy classes of
based maps
$$
 S^1\vee S^3 \;=\; S^1\times S^2/ \{0\}\times S^2\quad \longrightarrow
\quad \SO(3),
$$
i.e. to elements of 
$$
\pi_0(\Map_*(S^1\vee S^3, \SO(3))) = \Z/2\Z \oplus \Z,
$$
where $\Map_*$ is the based mapping space.  This has generators 
$$
\tau' \in \pi_0(\Map_*(S^1, \SO(3))) = \Z/2\Z,\qquad \al_0 \in
\pi_0(\Map_*(S^3, \SO(3))) = \Z.
$$
Thus $\pi_1(\Dd) \cong \Z\oplus\Z/2\Z \oplus \Z/2\Z$.
The two generators $\tau, \tau'$ of order $2$ lie in $\pi_1(K_0^0)$,
and it is not hard to check that $\al_0 = \al_1 + \tau'$, so that $\al_0\in
\pi_1(K_1^0)$  as claimed.

It follows that the composite
$$
G_\la^0\longrightarrow G_\infty^0\stackrel{\simeq}\longrightarrow \Dd
$$
induces an isomorphism on $\pi_1$.  This means that
 we can work out the relations
between the $\al_k$ in $\pi_1(\Dd)$, i.e. we do
not need to insist that all the diffeomorphisms we consider are symplectic. 
In particular we can alter $C_{2k}$ by an isotopy so that it has
the  form 
$$
C_{2k} = \{(z, \rho_k(z)): z\in S^2\},
$$
 where  $\rho_k: S^2\to S^2$
 is any convenient degree $-k$ map.  (This corresponds to conjugating
the original representative of $\al_{k}$ by an element of $\Dd$.)
Since the loop $\al_k$
 is represented by rotations of the fiber that fix
$C_{2k}$, we can represent $\al_k$ by  the loop
$$
\th\cdot (z,w) =  (z, R_{(\rho_k(z), \th)} w),\quad\th\in S^1 = \R/\Z,
$$
where $R_{(w,\th)}$ is the rotation by $\th$ about the axis through $w$.
It is easy to see that the corresponding element in
$\pi_0(\Map_*(S^1\vee S^3, \SO(3)))$ is $\tau' + k\al_0$, since the 
induced map $S^3\to \SO(3)$ has degree $k$. 
This implies that
$\al_k = k \al_1 + (k-1) \tau'$ as claimed.

To prove (ii), we must first consider the effect of the map
$G_\la^0\to \Dd$  on $H_3$. Now $H_3(\Dd) =
\Q\oplus \Q$ with generators $\xi_0, \eta$ from $ H_3(K_0^0)$. 
Since $H^2(U_0) = H^3(U_0) = 0$ by Corollary~\ref{cor:mv1},  
Proposition~\ref{prop:tensor} implies that
$H_3(G_\la^0) \cong H^3(K_0^0)$.  Hence the map $G_\la^0\to \Dd$ induces an
isomorphism on $H_3$.  Thus, again, we can work out the relations between
the $\xi_k$ in $\Dd$.

Since the $\SO(3)$-factor in $K_k^0$ is a
lifting to $S^2\times S^2$ of the usual action on $S^2$,
 $\xi_k$ is represented by a map of the form
$$
\Psi:\SO(3)\to G_\la:\quad  \Psi(g)\cdot(z,w) = (g(z), h_{k,g,z}w)),\;\; 
g, h_{k,g,z}\in
\SO(3). $$
Moreover, because all these diffeomorphisms map
the section 
$$
C_{2k} = \{(z, \rho_k(z)):z\in
S^2\}
$$
  into itself,
\begin{equation}\label{eq:rho}
h_{k,g,z}(\rho_k(z)) = \rho_k(g(z)).
\end{equation}
By our earlier remarks, $\rho_k$ can be any map of degree $-k$.  In particular we
may suppose that it has a fixed point $z_0$.  Let $e: \SO(3)\to S^2$ be the
evaluation map at $z_0$ and define $F_k$ by
$$ 
F_k: \SO(3)\to \SO(3):\quad g\mapsto h_{k,g, z_0}.
$$
Equation~(\ref{eq:rho}) implies that   there is a commutative diagram
$$
\begin{array}{ccc}  \SO(3) & \stackrel{F_k}\longrightarrow &
\SO(3)\\
e \downarrow & & e \downarrow\\
 S^2 &\stackrel{\rho_k}\longrightarrow & S^2
\end{array}
$$
Thus $F_k$ preserves the fibers of
$e$ and covers a map on the base of degree $-k$.  Since it has to commute
with the boundary map of the homotopy long exact sequence of $e$, it must
also have degree $-k$ on the fiber.  Therefore $F_k$ has degree $k^2$.

Now recall that, for any topological group $G$, the
sum of any two spherical cycles  $\al,\be$ in $H_k(G)$ is represented by the
product map
$$
S^k\to G:\quad g\mapsto \al(g)\, \be(g).
$$
Hence, since the $\xi_k$ are spherical, the class $\xi_k - \xi_0$
is represented by the map 
$$
\SO(3)\to G_\la:\quad  g\cdot(z,w) = (z, h_{k,g,z}(w)),\;\; 
h_{k,g,z},\,g\in S^3, $$
and hence by the map
$$
\SO(3)\to \Map(S^2,\SO(3)): \qquad g\mapsto (z\mapsto h_{k,g,z}).
$$
Similarly, $y$ is represented by
$$
\SO(3)\to\Map(S^2,\SO(3)): \qquad g\mapsto (z\mapsto g). 
$$
To check the identity $\xi_k - \xi_0 = k^2\eta$ we just have to see that 
$\xi_k - \xi_0$ and $k^2\eta$ induce the same map on rational $H_3$. 
Since the second loop space $\Om^2(\SO(3))$ has the same rational
homotopy type as $S^1$, the evaluation map 
$$
\Map(S^2, \SO(3)) \to \SO(3):\quad f\mapsto f(z_0)
$$
induces an isomorphism on $H_3(\cdot;\Q)$.  Thus the desired conclusion
follows from the fact that $F_k$ has degree $k^2$.
\QED

Next let us consider the analogous questions for the case $i = 1$. 
Here it is convenient to identify  $\cp\#\,\bcp$ with the set
$$
\{(z,w)\in D^2\times S^2 : (e^{i\th}, w) \equiv (1,R_{-\th} w), e^{i\th} \in \p
D^2\},
$$
where  $R_\th$ denotes the rotation of $S^2$ about its vertical axis
by  the angle $\th\in \R/2\pi\Z$, and we identify $D^2/\p D^2$ with $S^2$
so that the point $\p D^2$ corresponds to the south pole $P_s$. Then one can check
that the section $D^2\times \{P_s\}$ has self-intersection $-1$ so that it can be
identified with $C_1$.

 In this case
all groups $K_k^1$ are isomorphic to $\U(2)$. 
Further, the group $\Dd = \Dd^1$ of orientation preserving 
fiberwise diffeomorphisms of $\cp\#\,\bcp$ is homotopy
equivalent to the semi-direct product 
$$
\SU(2)\times_{\Z/2\Z} \Ga(S^2;\SO(3)),
$$
where the $\SU(2)$ factor is identified with a subgroup of $K_0^1$,
$\Ga(S^2;\SO(3))$ denotes the space of sections of the
bundle over $S^2$ with fiber $\SO(3)$ that is associated to 
the fibration $\cp\#\,\bcp\to S^2$, and $\Z/2\Z$ acts in a way that we now
explain.  If $R_\th$ denotes the rotation of $S^2$ about its vertical axis
by  the angle $\th\in \R/2\pi\Z$, it is easy to see that
$$
\Ga(S^2;\SO(3)) = \{\al\in \Map(D^2, \SO(3)): \al(e^{i\th}) = R_\th\al(1)
R_{-\th}, e^{i\th}\in \p D^2\}.
$$
With this identification, the generator of $\Z/2\Z$ acts by conjugation by
$R_\pi$ in $\Ga(S^2;\SO(3))$  and by multiplication by $-{\rm Id}$ in
$\SU(2) $.   To see this, compare  the present description
of $K_k$ with that given just before Proposition~\ref{prop:reduction}.
To make the two descriptions more alike, think of 
 $\Ga(S^2;\SO(3))$ as the quotient of $\Ga(S^2;\SU(2))$ by the action
of $\Z/2\Z$.  Then the nontrivial action of $\Z/2\Z$ on both factors of
$\SU(2)\times \Ga(S^2;\SO(3))$ described above
corresponds to the fact that the intersection of $\SU(2)\times S^1$ with $\T^2$
does not split, even after quotienting out by the $\Z/2\Z$ subgroup $\pm1$
of $S^1$.

 Denote by $\al_k, k\ge 0,$ the
generator of $\pi_1(K_k^1)$ corresponding to the $S^1$ action that fixes the
sphere $C_{2k+1}$, and by $\xi_k$ the generator of $H_3(K_k^1;\Z)$ coming
from the inclusion of $\SU(2)$ into $\U(2)$. Further, set $\eta = \xi_1 - \xi_0$.
\MS

\begin{lemma}\label{le:htpyK1}
\begin{itemize}
\item[{\rm (i)}] $\al_k = (k+1) \al_0 \in \pi_1(G_\la^1)$,

\item[{\rm (ii)}] $\xi_k = \xi_0 + k^2\eta \in H_3(G_\la^1; \Z)$.
\end{itemize}
\end{lemma}
\proof{}  The proof of (i) is left to the reader. 

Note that there is a fibration
$$
\Om^2(\SO(3)) \longrightarrow \Ga(S^2;\SO(3))\stackrel{ev}
 \longrightarrow \SO(3).
$$
where $ev(\al) = \al(1)$.  Hence $ev$ induces an isomorphism on $H_3$,
and so $H_3(\Dd)$ has again rank $2$. The formula in (ii) now follows
as in the previous lemma.  Observe that the analog of the set of elements
$h_{k,g,z_0}, g\in \SO(3),$ is a section of a nontrivial bundle.  However, if
we take $z_0$ to be the point $1\in \p D^2$ (which corresponds to the south
pole $P_s$ in $S^2$) and assume as before that the maps $\rho_k$ fix $z_0$, then the
composites  $h_{0,g,g(z_0)}^{-1}\circ h_{k,g, z_0}$ are elements of $\SU(2)$
and the argument can proceed as before.\QED

\section{Algebraic Computations}
\label{sec:algebra}	

This section is devoted to the algebraic computations necessary to
prove the main theorems of this paper. In \S\ref{ss:0} we analyse
 the case of $M^0_\la = (S^2 \times S^2, \om^0_\la)$ in complete detail.
The case $i = 1$ is  treated in
\S\ref{ss:1}.  It is completely analogous to the case $i = 0$.

\subsection{Computation of $H^{\ast}(G^0_\la)$}\label{ss:0}

As before, we assume
$\ell-1<\la\leq \ell$.   Since this whole subsection is devoted to $M^0_\la =
(S^2 \times S^2,  \om^0_\la)$, we will omit the superscript $^0$ to simplify
notation. Unless noted otherwise, we assume rational coefficients throughout.

Our strategy is first to show that there is a subalgebra of $H^*(G_\la)$
isomorphic to $\La(a,x,y)\otimes S(w_\ell)$ and then to show that 
this subalgebra must be the whole cohomology ring.  
Thus the first goal is to show how the new element $w_\ell$  appears.
To do this we use the
stratification of the space $\Jj_\la$ to construct an  interconnected family of
modified Mayer-Vietoris sequences $(MV_k)$   relating the cohomology
groups of the different strata.  This allows us to pinpoint exactly 
where $w_\ell$ comes from: see Definition~\ref{def:w}.  The
geometric significance of this generator is discussed further in
\S\ref{ss:kron}. In order to show that $\La(a,x,y)\otimes S(w_\ell)$ is the
whole of $H^*(G_\la)$, we  calculate all the maps in the part of 
$(MV_k)$ generated by this subalgebra, and then show that if $d$ is the
minimal degree in which $H^*(G_\la)$ contains something else, there must in
fact be a new element in some degree $< d$. 

Key to the whole process is the fact that each stratum $U_k$ is 
(weakly homotopic to)
a homogeneous space for $G_\la$.  Our first proposition shows that
the corresponding decompositions for $H^*(G_\la)$ are compatible with its
multiplicative as well as its additive structure.

\begin{prop}\label{prop:tensor}
\begin{itemize}
\item[{\rm (i)}] $H^{\ast}(G_\la)$ is a free algebra.  Moreover, for
all  $0\leq k\leq \ell,$
\begin{equation} \label{eq:tensor}
H^{\ast}(G_\la) \cong H^{\ast}(U_k) \otimes H^{\ast}(K_k),
\end{equation}
as graded algebras.

\item[{\rm (ii)}]  For all $k > 1$, $H^{\ast}(U_k)\cong H^{\ast}(U_1)$.  
\end{itemize}
\end{prop}
\proof{} Since $G_\la$ is a group, in
particular an $H$-space, we have that $H^{\ast}(G_{\lambda})$ is a Hopf
algebra. The Leray Structure Theorem for Hopf algebras over a field of
characteristic $0$ then says that $H^{\ast}(G_{\lambda})$ is the tensor
product of exterior algebras with generators of odd degree and polynomial
algebras with generators of even degree (see~\cite{Wh} for relevant
definitions and a proof of this structure theorem). Hence
$H^{\ast}(G_\la)$ is a free algebra.

The algebra isomorphism~(\ref{eq:tensor}) was proved in~\cite{Ab} for
$k=0$, and the proof for general $k$ is very similar.  The first point is that
the inclusion $i:K_k \rightarrow G_\la$ is injective in homology.
The result for $H_1(K_k)$ holds because, as in the proof of 
Lemma~\ref{le:htpyK0},
the generator $t_k$ of $\pi_1(K_k)\otimes \Q$ does not vanish in the space of
homotopy equivalences of $S^2\times S^2$.  The result for $H_3(K_k)$ is
somewhat subtle, and uses the fact that $G_\la$ maps onto
the space of homotopy equivalences of each sphere factor by maps of the
form:
$$
G_\la \to \Map(S^2, S^2):  \quad g\mapsto pr_1\circ g|_{S^2\times \{pt\}}.
$$
For more details, see~\cite{Ab}.

The Leray-Hirsch Theorem then implies that the Leray-Serre cohomology
spectral sequence for the fibration
$$
K_k\longrightarrow G_\la \stackrel{p}\longrightarrow G_{\lambda}/K_k
$$
collapses with
\[ E_{\infty}^{\ast,\ast}\cong E_{2}^{\ast,\ast}\cong
H^{\ast}(G_{\lambda}/ K_k)\otimes
H^{\ast}(K_k) \]
as bigraded algebras. This does not prove~(\ref{eq:tensor}),
but says that the bigraded algebra $E_{0}^{\ast,\ast}(H^{\ast}(G_\la),F)$,
associated to $H^{\ast}(G_{\lambda})$ with filtration $F$ coming from the
above fibration, is isomorphic to
$H^{\ast}(G_\la / K_k)\otimes H^{\ast}(K_k)$.

To finish the proof of~(\ref{eq:tensor}) we can now do the following.
$H^{\ast}(G_\la)$ has a subalgebra
$p^{\ast}(H^{\ast}(G_{\lambda}/ K_k))\cong
H^{\ast}(G_{\lambda}/ K_k )$.
Since
\[ i^{\ast}:H^{\ast}(G_{\lambda})\rightarrow H^{\ast}( K_k ) \]
is surjective, we can choose $a, b\in H^{odd}(G_\la)$ so that
$i^{\ast}(a)$ and $i^{\ast}(b)$ generate the ring
$H^{\ast}(K_k)$. The subalgebra of $H^{\ast}(G_{\lambda})$
generated by $a$ and $b$ is isomorphic to $H^{\ast}( K_k )$
since the only relations in this algebra hold universally
on any cohomology algebra (skew-symmetry of the cup product of elements
of odd degree).
Composing these inclusions of
$H^{\ast}(G_{\lambda}/ K_k)$
and $H^{\ast}(K_k)$ as subalgebras of
$H^{\ast}(G_{\lambda})$ with cup product multiplication in
$H^{\ast}(G_{\lambda})$ we get a map
\[ \nu:H^{\ast}(G_{\lambda}/ K_k)\otimes
H^{\ast}(K_k) \rightarrow H^{\ast}(G_{\lambda})\ .\]
Because $H^{\ast}(G_{\lambda})$ is (graded) commutative we have
that $\nu$ is an algebra homomorphism.
Moreover, by construction, $\nu$ is compatible with
filtrations (the obvious one on
$H^{\ast}(G_{\lambda}/ K_k)\otimes H^{\ast}( K_k )$
and $F$ on $H^{\ast}(G_{\lambda})$). As was already remarked, the
degeneration of the spectral sequence at the $E_{2}$ term implies that
the associated bigraded version of $\nu$ is an algebra isomorphism,
which in turn shows that $\nu$ itself is an algebra isomorphism,
concluding the proof of~(\ref{eq:tensor}).

Part (ii) of the proposition is a direct consequence of (i).\QED

\begin{cor}\label{cor:hlink} If $1\le k\le \ell-1$ the link $L_{k+1}$ of
$U_{k+1}$ in $U_k$ represents a nonzero element of $H_3(U_k)$.
\end{cor}
\proof{}  By Proposition~\ref{prop:link2}, $L_{k+1}$ is a homogeneous space 
for  $K_{k+1}$.  Hence, under the map $G_\la\to G_\la/K_k = U_k$,
the generator  of
$H_3(K_{k+1})$ is taken to the fundamental class of $H_3(L_{k+1})$.
Since the inclusion map $H_3(K_{k+1})\to H_3(G_\la)$ is nonzero, 
the same holds for the inclusion $H_3(L_{k+1}) \to U_k$.\QED

The computation of $H^{\ast}(G_\la)$
is based on using the isomorphisms in (i) above and
on the Mayer--Vietoris sequences for the decompositions
$$
U_{0\dots k} =  U_{0\dots \,k-1} \cup \Nn(U_k)
$$
where  $\Nn(U_k)$ is a neighborhood of $U_k$ in $U_{0\dots k}$.  By
Proposition~\ref{prop:link} we may choose $\Nn(U_k)$ to be a disc bundle over
$U_k$.  In particular, there is a fibration
$$
S^{4k-3}\to \p \Nn(U_k)\to U_k\ .
$$

\begin{lemma}\label{lem:zero}   For all $ k\leq \ell$
$$
H^i(U_{0\dots k})= 0,\qquad 0< i \leq 4k.
$$
\end{lemma}
\proof{}  The Mayer--Vietoris sequence for the decomposition
 $U_{0\dots k} = U_{0\dots\,k-1} \cup \Nn(U_k)$
is
$$
\to H^i(U_{0\dots k-1}) \oplus H^i(\Nn(U_k))\to
H^{i}(U_{0\dots k-1} \cap\Nn(U_k))\to H^{i+1}(U_{0\dots k})\to\cdots
$$
But, when $i < 4k-3$,
$$
H^{i}(U_{0\dots k-1} \cap\Nn(U_k)) = H^i(\p\Nn(U_k)) =
H^{i}(\Nn(U_k)),
$$
since the fiber of the projection $\p\Nn(U_k)\to U_k$ is connected up
to dimension $4k-3$.  This implies  that
$$
H^i(U_{0\dots k}) \cong H^i(U_{0\dots {k-1}}),\quad i \le 4k-4.
$$
The lemma follows by downwards induction on $k$ because  $U_{0\dots\ell}
= \Jj_\la$ is contractible.\QED

The next result is contained in~\cite{Ab}.  Its proof is repeated here for the
convenience of the reader.

 \begin{proposition}\label{lem:mv1}
\begin{itemize}\item[{\rm (i)}]
$H^1(U_0)=\Q$ and its generator is nonzero on the class in
$H_1(U_0)$ represented by the link $L_1$ of $U_{1}$ in $U_0$.

\item[{\rm (ii)}]  There is an exact  sequence
$$
(MV_1):\ \cdots\to 
H^{i-2}(U_1)\to H^i(U_{01})\to H^i(U_0)\stackrel{\tau}{\to}
H^{i-1}(U_1)
\stackrel{\de}{\to}
H^{i+1}(U_{01})\to\cdots
$$
\end{itemize}
\end{proposition}
\proof{}
Using Lemma~\ref{lem:zero} on the Mayer-Vietoris sequence for the
decomposition $U_{01} = U_0 \cup \Nn(U_1)$ we get the short
exact sequence:
$$
H^1 (U_{01}) = 0 \to H^1(U_0)\oplus H^1(\Nn(U_1)) \to
H^1(\p\Nn(U_1)) \to 0 = H^2(U_{01})
$$
If $r = \mbox{rk\ } H^1(\Nn(U_1)) = \mbox{rk\ } H^1(U_1)$,
it follows from Proposition~\ref{prop:tensor} (i) that
$r+1 = \mbox{rk\ } H^1(U_0)$ and hence we get from the previous
short exact sequence that $2r+1 = \mbox{rk\ } H^1(\p\Nn(U_1))$.
The existence of a fibration $S^1 \hookrightarrow \p\Nn(U_1) \to
U_1$ shows that
$$ 
2r+1 = \mbox{rk\ } H^1(\p\Nn(U_1)) \leq \mbox{rk\ } H^1(U_1) + 1
= r + 1 
$$
which implies that $r = 0$.

We conclude that $H^1(\p\Nn(U_1)) \cong \Q \cong  H^1(U_0)$  and
also that the inclusion $S^1 \hookrightarrow \p\Nn(U_1)$ is injective
in homology (or surjective in cohomology), which proves (i).

Moreover, it also follows that the cohomology spectral sequence for
the fibration
$$
S^1\to \p\Nn(U_1)\to U_1
$$
collapses and
$$
H^{\ast}(\p\Nn(U_1)) \cong H^{\ast}(S^1\times U_1)\cong
H^{\ast}(U_1)\oplus H^{\ast -1}(U_1).
$$
Part (ii) now follows easily.\QED

\begin{cor}\label{cor:mv1}
$H^2(U_0) = H^3(U_0) = 0$ and
$H^3(U_k) = \Q$ for $k\geq 1$.
Moreover, the smallest integer $d > 3$ such
that $H^d (U_k)\ne 0$ is the same for all $k\ge 0$.
\end{cor}
\proof{}  Since $H^1(U_0) = \Q$ we have by the isomorphisms in
Proposition~\ref{prop:tensor} (i) that $H^1(U_1) = 0$. Using this
and the fact that $H^i(U_{01}) = 0,\  0 < i\le 4$  by
Lemma~\ref{lem:zero}, in $(MV_1)$ we conclude that
$H^2(U_0) = H^2(U_1) = H^3(U_0) = 0$. From $H^3(U_0) = 0$,
and using again Proposition~\ref{prop:tensor} (i), we also get that
$H^3(U_k) = \Q$ for $k\geq 1$. The statement about
$d$ follows similarly by Proposition~\ref{prop:tensor} (i). \QED

\begin{proposition}\label{prop:mvk}
For all $1\le k \le \ell$ there is an exact sequence
$$
(MV_k):
\cdots\to H^i(U_{0\dots k})\to H^i(U_{0\dots
k-1})\stackrel{\tau}{\to} H^{i-4k+3}(U_k)\stackrel{\de}{\to}
H^{i+1}(U_{0\dots k})\to \cdots\ . $$
\end{proposition}
\proof{}
Observe that because $H^i(U_{0\dots k}) = 0$,
for $i\le 4k$, the Mayer--Vietoris sequence for the pair
$(U_{0\dots k}, U_{k})$ implies that
$$
H^i(U_{0\dots k-1}) \oplus H^i(\Nn(U_{k}))\cong
H^{i}(\p\Nn(U_{k})),\qquad i < 4k.
$$
But $\Nn(U_k) \simeq U_k$ and the fiber of the projection
$$
f: \p\Nn(U_{k})\to U_{k}
$$
is $S^{4k-3}$.
Therefore the fact that $ H^{4k-2}(U_k) $ injects into
$ H^{4k-2}(\p\Nn(U_k)) $ implies that
all differentials in the Leray spectral sequence 
of $f$ vanish.  Hence the map
$$
H^i(U_k) \cong  H^i(\Nn(U_k)) \to H^i(\p\Nn(U_k))
$$
is injective for all $i$ (not just for $i < 4k$) and has cokernel isomorphic to
$H^{i-4k+3}(U_k)$.  Hence $(MV_k)$ is exact for all $k\le \ell$. \QED

\begin{proposition}\label{prop:deltak}
\begin{itemize}\item[{\rm (i)}] In $(MV_k)$, for $k < \ell$, the boundary map
$$
\de_k: \Q = H^3(U_k)\to H^{4k+1}(U_{0\dots k})
$$
is an isomorphism.

\item[{\rm (ii)}]  $H^{4k+j}(U_{0\dots k}) = 0$ for
$j = 2,3$ and $1\le k < \ell $. Moreover
\begin{eqnarray*}
H^{4\ell}(U_{0\dots \ell -1}) & = & \Q,\\
H^{4k+4}(U_{0\dots k}) & = & 0,\quad k < \ell - 1.
\end{eqnarray*}
\end{itemize}
\end{proposition}
\proof{} (i)
>From $(MV_{k+1})$ with $i = 4k+1$ we find that
$$
H^{4k+1}(U_{0\dots {k+1}})\to H^{4k+1}(U_{0\dots\, {k}})
\stackrel{\tau}{\to} H^0(U_{k+1}) \to H^{4k+2}(U_{0\dots {k+1}})
$$
is exact.  Since the first and last terms above are $0$ by Lemma~\ref{lem:zero}
and the third is $\Q$, we get that $ H^{4k+1}(U_{0\dots {k}})= \Q$.
Moreover, the way $(MV_{k+1})$ is constructed out of Mayer-Vietoris for
$U_{0\dots\, {k+1}} = U_{0\dots\, {k}} \cup \Nn(U_{k+1})$ and degeneration
of
$$
S^{4k+1}\to \p\Nn(U_{k+1})\to U_{k+1}
$$
shows that the generator of $H^{4k+1}(U_{0\dots {k}})$ is non-zero on the
link $S^{4k+1}$ of $U_{k+1}$ in $U_{0\dots {k}}$, and so this link is 
an explicit geometric generator for the homology group 
$H_{4k+1}(U_{0\dots {k}})$.

The homology version of the map $\de_k: H^3(U_k)\to H^{4k+1}(U_{0\dots k})$
in $(MV_{k})$ sends the link $S^{4k+1}$ of $U_{k+1}$ in $U_{0\dots {k}}$
to its intersection with $U_{k}$, i.e. to the link of $U_{k+1}$ in
$U_{k}$. Hence,  $\de_{k}$ is
an isomorphism because, by Corollary~\ref{cor:hlink}, the link of $U_{k+1}$ in
$U_{k}$ represents a non-zero element in $H_{3}(U_{k})$. 

To prove (ii), look at $(MV_{k+1})$ with $i = 4k + j$, where $j=2,3,4$, using
the facts that $H^1(U_k) = H^2(U_k) = 0$ and that
$\de_{k+1}$ is an isomorphism when $k+1 < \ell$. 
 \QED

Our next aim is to  show that the nonzero element in
$H^{4\ell}(U_{0\dots\ell-1})$ restricts to give a nonzero element in
$H^{4\ell}(U_0)$.  Moreover, we will see that this  is the ``first" nonzero
element in the $H^*(U_k)$ in degrees $* > 3$.  More precisely, as in
Corollary~\ref{cor:mv1}, we consider the smallest integer $d = d_\ell > 3$
 such that $H^d(U_k)\ne 0$ for any $k$, and will show that $d_\ell =
4\ell$.

\begin{lemma}\label{lem:d2} $d_1 = 4$, $d_2 = 8$ while  $d_\ell > 11$
when  $\ell > 2$.
\end{lemma}
\proof{} The fact that $d_1 = 4$ follows from Proposition~\ref{prop:deltak}
(ii).

Consider the sequence $(MV_1)$ with $i = d_\ell = d$:
$$
\cdots\to
H^{d-2}(U_1)\stackrel{\de}{\to} H^d(U_{01})\to
H^d(U_0)\stackrel{\tau}{\to}  H^{d-1}(U_1) = 0\ .
$$
Clearly, $d$ is the smallest integer $> 3$ such that
$H^d(U_{01})$ is nonzero and also not covered by the image of
$\de$.

Now suppose that $\ell = 2$.  In this case the sequence $(MV_2)$
collapses to the isomorphisms
$$
\tau: H^i(U_{01})\stackrel{\cong}{\to} H^{i-5}(U_{2}), i>0.
$$
Hence the first nonzero groups $H^i(U_{01})$ (for $i > 0$) occur
when $i = 5,8$. By Proposition~\ref{prop:deltak} (i) we know that
$H^5(U_{01})$ is in the image of $\de$, and so $d_2 \ge 8$. This
implies that $H^6(U_1) = 0$, which means that $H^8(U_{01})$ is
not in the image of $\de$.  Hence $d_2 = 8$.

If $\ell > 2$ the argument showing that $d_\ell > 5$ goes through
as before. Now look at $(MV_2)$:
$$
\cdots\to H^i(U_{012})\to H^i(U_{01})
\stackrel{\tau}{\to}  H^{i-5}(U_2)
\stackrel{\de}{\to} H^{i+1}(U_{012})\to\cdots\ . $$
We know, by Lemma~\ref{lem:zero}, that $H^i(U_{012}) = 0$ for
$i\le 8$. Hence $H^i(U_{01}) = 0$ for these $i$, since in this
range $H^{i-5}(U_2)$ is either zero or mapped isomorphically by
$\de$. Hence $d_\ell > 8$. By Proposition~\ref{prop:deltak}(i) we
know that $ H^{9}(U_{012})= \Q$ is in the image of $\de$, while
Proposition~\ref{prop:deltak}(ii) says that   $ H^{i}(U_{012})=
0$ for $i = 10, 11$ (and also $12$ when $\ell > 3$).  But in this
range, because $H^{i-5}(U_2) = 0$, we get from $(MV_2)$ that
$H^{i}(U_{01})\cong H^{i}(U_{012})=0$.  Feeding this
information back into $(MV_1)$, we find that $d_\ell > 11$ as
claimed.  In fact, if $\ell > 3$ we get that $d_\ell > 12$. \QED

More generally:

\begin{prop}\label{prop:dl}  $d_\ell = 4\ell$, for any $\ell\in\N$.
\end{prop}
\proof{} We will prove by induction on $\ell\in\N$ the statement:

\begin{center}
for any $\ell\in\N$, $d_\ell = 4\ell$ and $d_p > 4\ell + 3$ for $p
> \ell$
\end{center}

\NI which implies the result and is more convenient for the inductive
step.

As the previous lemma shows, the statement is true for $\ell=1,2$.
We now suppose, by induction, that the statement is true for $\ell
- 1 \ge 2$ and will consider the case $\ell$.
Combining Lemma~\ref{lem:zero} with Proposition~\ref{prop:deltak}
we have that
\begin{eqnarray*}
H^i(U_{0\dots\,\ell-1}) & = & 0, \quad i < 4\ell, i\ne 4\ell-3;\\
H^{4\ell -3}(U_{0\dots\,\ell-1}) & \cong & \Q\\
H^{4\ell}(U_{0\dots\,\ell-1}) & \cong & \Q.
\end{eqnarray*}
We will show that for all $0 \le k < \ell -1$,
\begin{eqnarray*}
H^i(U_{0\dots\,k}) & = & 0, \quad 4\ell -4 < i < 4\ell\\
H^{4\ell}(U_{0\dots k}) & \cong & \Q.
\end{eqnarray*}
This for $k=0$ together with the inductive hypothesis saying that
$d_\ell > 4\ell - 1$ immediately implies that $d_\ell = 4\ell$.

Consider the sequence $(MV_k)$ for $k = \ell-1,
\ell-2$ and so on:
$$
\cdots\to H^i(U_{0\dots k})\to H^i(U_{0\dots
k-1})\stackrel{\tau}{\to} H^{i-4k+3}(U_k)\stackrel{\de}{\to}
H^{i+1}(U_{0\dots k})\to\cdots
$$
When $4\ell - 4 < i \le 4\ell$ and $1\le k \le \ell - 1$, the
numbers $i - 4k + 3$ range in the interval $[4, 4\ell -1]$, and so
the inductive hypothesis implies that all $H^{i-4k+3}(U_k)
= 0$. In order to deduce that
$H^i(U_{0\dots k})\cong H^i(U_{0\dots k-1})$, for $4\ell - 4 <
i \le 4\ell$ and $1\le k \le \ell - 1$, we also need to know that
$H^{i-4k+2}(U_k) = 0$. This holds except when $i = 4\ell - 3$ and
$k = \ell -1$, since then $i - 4k + 2 = 3$. However, this nonzero
group exactly cancels out the nonzero group $H^{4\ell-3}(U_{0\dots
\ell-1})$ via $\de$.  Hence the groups $H^i(U_{0\dots k}) = 0$ for
$ 4\ell-4 < i < 4\ell$  and $k = \ell -2$, and remain $0$ as $k$
decreases to $0$. Similarly,  the groups $H^{4\ell}(U_{0\dots k})
= \Q$ for all $0\le k \le \ell -1$.

Now consider the case $p-1 < \la \le p$ with $\ell < p \in\N$. In
this case we know that
\begin{eqnarray*}
H^i(U_{0\dots\,\ell-1}) & = & 0, \quad i \le 4\ell, i\ne
4\ell-3;\\ H^{4\ell -3}(U_{0\dots\,\ell-1}) & = & \Q,
\end{eqnarray*}
and the previous argument shows that $d_p > 4\ell$. We also must
show that $H^i(U_0) = 0$ for $4\ell < i \le 4\ell + 3$. But
because $p > \ell$ Proposition~\ref{prop:deltak} says that
\begin{eqnarray*}
H^i(U_{0\dots\ell}) & = & 0, \quad  i = 4\ell + 2,  4\ell + 3;\\
H^{4\ell + 1}(U_{0\dots\,\ell}) & = & \Q.
\end{eqnarray*}
If we repeat the previous argument (replacing $\ell$ by $\ell + 1$)
we find that $H^i(U_{0\dots k}) = 0$ for $4\ell + 1 \le i \le
4\ell + 3$ and $k = \ell-1, \ell -2,\dots, 1$.  We cannot use this
argument to reach $k = 0$.  However, as at  the beginning of the
proof of Lemma~\ref{lem:d2}, it is enough to know that
$H^i(U_{01}) = 0$ for these $i$ to conclude that  $d_p > 4\ell +
3$. \QED

\begin{defn}\label{def:w}\rm  We define an explicit  new element
$w_\ell\in H^{4\ell}(G_\la)$. Let $y\in H^*(G_\la)$ come from the
second $\SO(3)$ factor in $K_0$ as in Lemma~\ref{le:htpyK0}, and
let $y_k\in
H^3(U_k)$ be its image under  the isomorphism of
Proposition~\ref{prop:tensor}. 
Because $U_{0\cdots\ell}\equiv \Jj_{\la}$ is contractible, 
 the map 
$$
\tau: H^{4\ell}(U_{0\dots\ell-1}) \longrightarrow H^3(U_\ell)
$$
in $(MV_{\ell})$ is an isomorphism.  Hence  $y_\ell$ 
lifts to an element ${\ov{y}}_\ell$ in $H^{4\ell}(U_{0\dots\ell-1})$ and then,
by the arguments above, restricts to a nonzero element $v_0$ of  $H^{4\ell}(U_0)$. 
By Proposition~\ref{prop:tensor}, there is a corresponding element in
$H^{4\ell}(G_\la)$ and it is this that is called $w_\ell$.  When there is no
possibility of confusion, we will write $w$ instead of $w_\ell$.
Finally, for $k=1,\ldots,\ell$, we denote by $v_k$ the image of $w$ in
$H^{4\ell}(U_k)$ 
under the isomorphisms of Proposition~\ref{prop:tensor}.  \end{defn}

\begin{remark}\label{rmk:w}\rm
It is perhaps easier to see ${\ov{y}}_\ell$ in homology.  Let $\eta$ be
a spherical representative of the second $\SO(3)$ factor in $H_3(K_0)$, and let
 $S \cong S^{4\ell-3}$ be the 
 link of $U_\ell$ in $\Jj_\la = U_{0\dots\ell}$.  Because $K_0$ acts on $\Jj_\la$
by multiplication on the left, there is a well defined $4\ell$-cycle $\eta* S$
in $U_{0\dots \ell-1}$ and by definiton ${\ov{y}}_\ell(\eta*S) = 1$. (For more
details of this notation, see Lemma~\ref{le:sa2}.) Note that we can replace $\eta$
here by any $3$-cycle in $K_0$ that is nonzero  in $K_\ell$, i.e. by anything
except a multiple of $\xi + \ell^2\eta$.  In \S\ref{ss:kron}, it is natural to 
consider the cycle $\xi + (\ell + 1)^2\eta$ that lives in $K_{\ell+1}$, while in
\S\ref{ss:whit} we will consider $\xi$.
\end{remark}

Since $H^\ast (G_\la)$ is the product of an exterior algebra with a
polynomial algebra, we now know that
$$ 
\La(a,x,y)\otimes S(w) \subset H^\ast (G_\la),
$$
with $\mbox{deg}\,a=1,\ \mbox{deg}\,x = \mbox{deg}\,y = 3 $ and 
$\mbox{deg}\,w = 4\ell$. 
This completes the first part of the argument.
Our final task is to show that these are all
the generators of $H^\ast(G_\la)$, i.e., the above inclusion is 
actually an equality of algebras.
To do this, we first  calculate all the groups and maps that occur in  the part
of  $(MV_k)$ generated by this subalgebra, showing in particular that the
relevant part of $H^*(U_{0\dots k})$ is a free algebra: see
Lemma~\ref{le:free}.  Then, we
 show that if $d$ is the
minimal degree in which $H^*(G_\la)$ contains something else, there must in
fact be a new element in some degree $< d$.

First we must work out the relations between the different elements
$v_k\in H^{4\ell}(U_k)$ that were defined in Definition~\ref{def:w}.
We just saw in the proof of the 
previous proposition that the restriction maps
$$
 H^i(U_{0\cdots k}) \to H^i(U_{0\cdots k-1}) 
$$
are isomorphisms when $1\leq k \leq \ell-1$ and $4\ell-4\leq i\leq
4\ell$. Hence
$$ 
i^\ast : H^{4\ell}(U_{0\cdots k}) \to H^{4\ell}(U_0) 
$$
is an isomorphism for $0\leq k\leq \ell-1$, and we will denote by
$v_0$ the generator of each $H^{4\ell}(U_{0\cdots k})$ such that
$i^\ast(v_0)=v_0$.

\begin{lemma} \label{lem:ik} Let $i_k$ denote the inclusion
$U_k \to U_{0\cdots k}$. Then  $i_k^\ast (v_0) = v_k$, for $0\leq k\leq
\ell-1$. \end{lemma}
\proof{}  This holds for $k = 0$ by definition.  
To avoid confusion in the following argument, we will think of 
 $v_0$ as an element of
$H^*(U_{0\dots k})$ and will write $i_0^*(v_0)$ for its restriction to $U_0$.
Suppose $k > 0$, and
let $f:P\to G_\la$ be a $4\ell$-dimensional cycle in
$G_\la$ on which $w$ does not vanish. It follows from the construction
of the isomorphisms in Proposition~\ref{prop:tensor} that the cycle
$f_k:P\to U_k$, defined by $f_k(p)=f(p)_\ast (J_{2k})$ where $J_{2k}$ is the
canonical complex structure in $U_k$, satisfies $v_k(f_k)=w(f)\ne 0$, for
all $1\leq k\leq \ell$.  Further, $i_0^*(v_0)(f_0) = w(f)$.

Since there is a path of almost complex structures in $U_{0\cdots k}$
joining $J_{2k}$ to $J_0$, the cycles $f_k$ and $f_0$ are
homotopic in $U_{0\cdots k}$. Hence, 
 $$
v_0(f_k)=v_0(f_0)= i_0^*(v_0)(f_0) = w(f)=v_k(f_k).
$$
Since $H^{4\ell}(U_k)\cong\Q$, this implies that $i_k^\ast(v_0)=v_k$
as claimed. \QED

The collapse of the spectral sequence of the fibration
$$ 
S^{4k-3}\to\p\Nn(U_k)\stackrel{\pi}{\to} U_k,
$$
together with the fact that the fiber is an odd-dimensional sphere,
implies that we have ring isomorphisms
$$
H^{\ast}(\p\Nn(U_k))\cong H^{\ast}(S^{4k-3})
\otimes H^{\ast}(U_k),\ 1\leq k\leq\ell.
$$
Hence we know that
$$
H^{\ast}(\p\Nn(U_k))\supseteq \La(\pi^\ast(y_k),
\widehat{z}_{4k-3})\otimes S(\pi^\ast(v_k)),\ 1\leq k\leq\ell,
$$
where
\begin{itemize}
\item[{\rm (i)}] $y_k \in H^{3}(U_k)$ is as  in
Definition~\ref{def:w};
\item[{\rm (ii)}] $\widehat{z}_{4k-3}\in H^{4k-3}(\p\Nn(U_k))$
is the restriction of $z_{4k-3}\in H^{4k-3}(U_{0\cdots k-1})$, and
$z_{4k-3}$ is determined from the isomorphism
$$
\tau: H^{4k-3}(U_{0\cdots k-1})\to H^{0}(U_k)
$$
 in $(MV_k)$ by
the condition $\tau(z_{4k-3}) = 1$.
\end{itemize}

\begin{lemma} \label{lem:rest} Let $j_k$ denote the inclusion 
$\p\Nn(U_k) \to U_{0\cdots k-1}$ and $\pi$
the projection $\p\Nn(U_k)\to U_k$. Then:
\begin{itemize}
\item[{\rm (i)}] $j_k^{\ast}(v_0) = \pi^{\ast}(v_k)$, for $0\leq k\leq\ell-1$;
\item[{\rm (ii)}] $j_\ell^{\ast}(v_0) = \widehat{z}_{4\ell - 3}\cdot
\pi^{\ast}(y_\ell) + c \pi^{\ast}(v_\ell)$, with $c\neq 0$.
\end{itemize}
\end{lemma}
\proof{} Statement (i) follows trivially from the previous lemma. Since $v_0$
lives in $U_{0\cdots k}$, its restriction to $\p\Nn(U_k)$ is
given by $\pi^\ast(i_k^\ast v_0) = \pi^\ast(v_k)$.

To prove (ii), observe first that the isomorphism $\tau: H^{4\ell}(U_{0\cdots
\ell-1}) \to H^{3}(U_\ell)$ from $(MV_\ell)$ is the composite of the
restriction map 
$$
j_\ell^\ast : H^\ast(U_{0\cdots \ell-1})\to H^\ast
(\p\Nn(U_\ell))$$
with integration over the fiber of the projection
$$\pi :\p\Nn(U_\ell)\to U_\ell\ .$$
Since integration over the fiber kills the elements of $\pi^\ast
(H^\ast(U_\ell))$ and $\tau(v_0)=y_\ell$, we must have
$$j_\ell^\ast (v_0) = \widehat{z}_{4\ell-3}\cdot \pi^\ast(y_\ell) + 
c \pi^\ast(v_\ell)\ .$$
If $c=0$ then $j_\ell^\ast(v_0)^2 = 
(\widehat{z}_{4\ell-3}\cdot \pi^\ast(y_\ell))^2 = 0$
because of degree considerations. Since $j_\ell^\ast$ is injective (because
$\tau$ is an isomorphism) this implies $(v_0)^2 =0$. But this is impossible,
since the restriction of $v_0$ to $H^\ast(U_0)$ is non-zero and
$H^\ast(U_0)$ is a free algebra. Hence $c\neq 0$ as stated.
\QED

\begin{lemma} \label{le:free} For each $k, 1\leq k\leq\ell$, the elements
$z_{4k-3}$ and $v_0$ generate a free subalgebra of 
$H^{\ast}(U_{0\cdots k-1})$.
\end{lemma}
\proof{} If $1\leq k\leq\ell-1$ this follows since
$$
j_k^{\ast}(z_{4k-3}\cdot v_0^n) = \widehat{z}_{4k-3}
\cdot (\pi^{\ast}(v_k))^n \neq 0 
$$
in $H^{\ast}(\p\Nn(U_k))$. For $k=\ell$ we have
$$j_\ell^\ast(z_{4\ell-3}\cdot v_0^n) = 
\widehat{z}_{4\ell-3}\cdot[\widehat{z}_{4\ell-3}\cdot \pi^\ast(y_\ell) + 
c \pi^\ast(v_\ell)]^n = c^n\cdot\widehat{z}_{4\ell-3}\cdot
[\pi^{\ast}(v_\ell)]^n \neq 0$$
in $H^{\ast}(\p\Nn(U_\ell))$.
\QED

\begin{proposition} \label{prop:maps}  Let $1\leq k\leq\ell-1$.
The maps $i^{\ast}$, $\tau$ and $\de$ in the exact sequence $(MV_{k})$
$$
\cdots\to H^i(U_{0\cdots k})\stackrel{i^{\ast}}{\to} H^i(U_{0\cdots
k-1})\stackrel{\tau}{\to} H^{i-4k+3}(U_k)\stackrel{\de}{\to}
H^{i+1}(U_{0\cdots k})\to\cdots
$$
satisfy:
$$
\begin{array}{lll}
i^{\ast}(v_{0}^n) = v_{0}^n & \tau(v_{0}^n) = 0 & \de(v_{k}^n) = 0 \\
i^{\ast}(z_{4k+1}\cdot v_{0}^n) = 0 & \tau(z_{4k-3}\cdot v_{0}^n)
= v_{k}^n & \de(y_{k}\cdot v_{k}^n) = c\cdot z_{4k+1}\cdot v_{0}^n,\ 
c\neq 0\ .
\end{array}
$$
\end{proposition}
\proof{} 

\begin{itemize}
\item[$i^{\ast}$] is induced by inclusion, hence a ring homomorphism. By
construction we know that $i^{\ast}(v_0) = v_0$. Because 
$\de:H^3(U_k)\to H^{4k+1}(U_{0\cdots k})$ is an isomorphism for
$k<\ell$ (by Proposition~\ref{prop:deltak} (i)), we also know that
$i^{\ast}(z_{4k+1})=0$.
\item[$\tau$] Since $i^{\ast}(v_0^{n}) = v_0^n$ we have that $\tau(v_0^n)=0$.
By construction $\tau(z_{4k-3})=1$. Also
$$\tau(z_{4k-3}\cdot v_0^n) = \int_{S^{4k-3}} \widehat{z}_{4k-3}
\cdot (\pi^{\ast}(v_k))^n = v_k^n\ .$$
\item[$\de$] Since $\tau(z_{4k-3}\cdot v_0^n) = v_k^n$ we have that
$\de(v_k^n)=0$. As before, since $\de:H^3(U_k)\to H^{4k+1}(U_{0\cdots k})$ 
is an isomorphism and $H^3(U_k)\cong  \Q$,
we know that $\de(y_k) = c\cdot z_{4k+1}$, with $c\neq 0$.
\end{itemize}

It remains to calculate $\de(y_{k}\cdot v_{k}^n)$.  This is where the
hypothesis $k\le \ell-1$ comes in.  We will restrict  $(MV_{k})$ to
a neighborhood of $U_{k+1}$ and exploit the fact that the
resulting sequence is a module over $H^*(U_{k+1})$.   The construction
uses the fact that  the $U_k$ form a stratification of $\Jj_\la$ with nice
normal structure: see the discussion before Proposition~\ref{prop:link}.
The neighborhoods occuring below are chosen to be fibered, with fiber equal
to a cone over the relevant link.

Recall that $(MV_k)$ is derived
from the  Mayer--Vietoris sequence for the decomposition $U_{0\cdots k} =
U_{0\cdots k-1}\cup \Nn(U_k)$.   Intersecting all sets with 
$\Nn (U_{k+1})$, we may identify $U_{0\cdots k}\cap \Nn (U_{k+1})$ with
its boundary $\p\Nn(U_{k+1})$,  and $\Nn (U_{k})\cap \Nn (U_{k+1})$ with
its deformation retract  $\p\Nn_{k}(U_{k+1})$.  (Here we write $\Nn_k(X)$
for the neighborhood of $X$ in $U_k$ and $\p\Nn_k(X)$ for its boundary.)  
Note also that 
$$
Y = U_{0\cdots k-1}\cap \Nn(U_k)\cap \Nn (U_{k+1})
$$
 is
(rationally) homotopy equivalent to an odd-dimensional sphere bundle over
$\p\Nn_{k}(U_{k+1})$ ($S^{1}$-bundle if $k=1$ and  $S^{3}$-bundle if $k>1$). 

Thus  the  Mayer--Vietoris sequence has the form 
$$ 
\to H^i(\p\Nn_{k}(U_{k+1}))\oplus H^i(U_{0\cdots k-1}\cap 
\Nn (U_{k+1}))\to H^i(Y)\to H^{i+1}(\p\Nn (U_{k+1}))\to 
$$
and, as with $(MV_k)$, it gives rise to an exact sequence
$$
(*)\ \to
H^{i}(U_{0\cdots k-1}\cap \Nn (U_{k+1})) \stackrel{\tau'}{\to} 
H^{i-4k+3}(\p\Nn_{k}(U_{k+1})) \stackrel{\de'}{\to}
H^{i+1}(\p\Nn(U_{k+1}))\to
$$
All these spaces fiber over $U_{k+1}$ in a compatible way.
Hence the homology groups involved are $H^{\ast}(U_{k+1})$-modules,
all the maps in the Mayer-Vietoris sequence preserve this module
structure, and it is not hard to check that those in
the modified sequence $(*)$ do as well.

Next, consider the commutative diagram:
$$ 
\begin{array} {ccc}
H^{i-4k-3}(U_{k}) &\stackrel{\de}{\to} & H^{i+1}(U_{0\cdots k})\\
\hr^*\downarrow & & \downarrow  r^*\\
H^{i-4k-3}(\p\Nn_{k}(U_{k+1})) &\stackrel{\de'}{\to} & H^{i+1}(\p\Nn 
(U_{k+1}))
\end{array}
$$
where  the  vertical maps are the obvious restrictions.  Thus $r^*$ is the map
previously called $j_{k+1}^*$.  From the two preceding 
lemmas we know that $r^*$ is injective and 
$$ 
r^* (\de(y_{k})) = c r^* (z_{4k+1}) =
c \widehat{z}_{4k+1},\ c\neq 0.
$$
Defining $\widehat{y}_{k} = \hr^*(y_{k})$, and using
commutativity of the diagram, we have that $\de'(\widehat{y}_{k})) =
c \widehat{z}_{4k+1}$.

In order to compute $\de'\circ \hr^*(y_{k}\cdot v_{k}^n)$
we will show that $\hr^*(v_{k})$ is pulled back from $H^*(U_{k+1})$
via the projection
$$
\pi: \p\Nn_k(U_{k+1}) \longrightarrow U_{k+1},
$$
 and
then will use the module property. Consider the commutative diagram of
inclusions: $$
\begin{array} {ccc}
\p\Nn (U_{k+1}) &\stackrel{r}{\to} & U_{0\cdots k}\\
i_\p \uparrow & & \uparrow i_{k} \\
\p\Nn_{k}(U_{k+1}) &\stackrel{\hr}{\to} & U_{k}
\end{array}
$$
>From Lemmas~\ref{lem:ik} and~\ref{lem:rest} we have that for
$1\leq k\leq \ell -2$
$$
\hr^{\ast}(v_{k}) = \hr^{\ast}\circ i_{k}^{\ast} (v_{0})
= i_\p^{\ast}\circ r^{\ast}(v_{0})
= i_\p^{\ast}\circ \pi^{\ast}(v_{k+1}).
$$
(Recall that $r = j_{k+1}$.)  On the other hand, for $k=\ell - 1$
$$
\hr^{\ast}(v_{\ell-1}) =
i_\p^{\ast}\circ r^{\ast}(v_{0}) =
i_\p^{\ast} (\widehat{z}_{4\ell-3}\cdot\pi^{\ast}(y_{\ell})
+ c \pi^{\ast}(v_{\ell})),\ c\neq 0.
$$
If we define
$$
\hpi=\pi\circ i_\p:\p\Nn_{k}(U_{k+1})\to 
U_{k+1}\ ,$$
this can be written as
$$
\hr^{\ast}(v_{k}) = \hpi^{\ast}(v_{k+1})\ 
\mbox{for}\ 1\leq k\leq \ell -2\ ,$$
while
$$
\hr^{\ast}(v_{\ell-1}) = i_\p^{\ast} 
(\widehat{z}_{4\ell-3})\cdot \hpi^{\ast}(y_{\ell}) +
c \cdot \hpi^{\ast}(v_{\ell}).
$$
We saw in Proposition~\ref{prop:link} that the fiber $L_\ell$ of the fibration
$$
L_{\ell}\hookrightarrow \p\Nn_{\ell-1}(U_{\ell}) \stackrel
{\hpi}{\to} 
U_{\ell}
$$
is a $3$-dimensional lens space.  Moreover the Leray--Serre spectral
sequence collapses by Corollary~\ref{cor:hlink}.  Hence 
$$
H^{4\ell}(\p\Nn_{\ell-1}(U_{\ell}))\cong 
\hpi^*(H^{4\ell}(U_{\ell}) \oplus H^{4\ell-3}(U_{\ell}) \cong
H^{4\ell}(U_{\ell}
$$
by Proposition~\ref{prop:dl}.  Thus
$$
\hr^{\ast}(v_{\ell-1}) = 
c \cdot \hpi^{\ast}(v_{\ell}).
$$

We can now use the fact that $\de'$ is a map of 
$H^{\ast}(U_{k+1})$-modules to conclude that
for $1\leq k\leq \ell -2$
\begin{eqnarray*}
\de'\circ \hr^{\ast}(y_{k}\cdot v_{k}^n) & = & 
\de'(\widehat{y}_{k}(\hr^{\ast}(v_{k}))^{n}) \  = \ 
\de'(\widehat{y}_{k} \hpi^{\ast}(v_{k+1})^{n}) \\
 & = & \pi^{\ast}(v_{k+1})^{n}\cdot \de'(\widehat{y}_{k})\  = \ 
r^{\ast}(v_{0})^{n}\cdot c \cdot \widehat{z}_{4k+1} \\
 & = & r^{\ast}(c \cdot z_{4k+1}\cdot v_{0}^{n}),
\end{eqnarray*}
while for $k=\ell -1$
$$
 \de'\circ \hr^{\ast}(y_{\ell-1}\cdot v_{\ell-1}^n) \ = \ 
r^{\ast}(c^n c\cdot z_{4\ell -3}\cdot v_{0}^{n})\ .$$
Therefore,
$$ 
\de(y_{k}\cdot v_{k}^n) = c\cdot z_{4k+1}\cdot v_{0}^n,\quad 
c\neq 0 $$
 for all $1 \le k \le \ell-1$, as claimed.
\QED

\NI
{\bf Proof of Theorem~\ref{thm:main0}}

Suppose that 
$$
H^{\ast}(G_{\la}) \ne A = \La (a,x,y) \otimes S(w),
$$
 and let
$$
u\in H^{\ast}(G_{\la}) - A
$$
be a non-zero element of minimal degree $d$.  From the arguments
leading to Proposition~\ref{prop:dl} and its proof we know that
$d>4\ell$.  The existence of this new element $u$ gives rise to new
non-zero elements
$$
u_k \in H^d (U_{0\cdots k}) - (\La(z_{4k+1}) \otimes S(v_0)),$$
for $0\leq k\leq \ell -1$, in the following inductive way. For
$k=0$ we let $u_0$ be the image of $u$ under the relevant
isomorphism from Proposition~\ref{prop:tensor} (i). Given a
non-zero element
$$u_{k-1} \in H^d (U_{0\cdots k-1}) -  (\La(z_{4k-3}) \otimes S(v_0)),
\ 1\leq k \leq\ell -1,$$
consider the following piece of $(MV_k)$:
$$\cdots\to H^d(U_{0\cdots k})\stackrel{i^{\ast}}{\to} H^d(U_{0\cdots
k-1})\stackrel{\tau}{\to} H^{d-4k+3}(U_k)\stackrel{\de}{\to}
H^{d+1}(U_{0\cdots k})\to\cdots
$$ 
By exactness we know that that
$\tau(u_{k-1})\in\,\mbox{Ker}\,\de$, and by
Proposition~\ref{prop:maps} and minimality of $d$ this implies that
$$\tau(U_{k-1}) = v_{k}^{n} = \tau (z_{4k-3}\cdot v_{0}^{n})$$
for some $n$. Hence the non-zero element $u_{k-1}-z_{4k-3}\cdot
v_{0}^{n}$ belongs to $\mbox{Ker}\,\tau$ and, using exactness
again, there must exist a non-zero $u_{k}\in H^{d}(U_{0\cdots k})$
such that
$$ 
i^{\ast}(u_k) = u_{k-1}-z_{4k-3}\cdot v_{0}^{n}.
$$
It is clear from Proposition~\ref{prop:maps} that $u_k
\not\in\La(z_{4k+1})\otimes S(v_0)$.

The element $u_{\ell -1}\in H^{d}(U_{0\cdots\ell-1}) - 
(\La(z_{4\ell -3})\otimes S(v_0))$ determines, via $(MV_\ell)$ and
the fact that $U_{0\cdots\ell}$ is contractible, a new non-zero
element
$$
u_{\ell}\in H^{d-4\ell +3}(U_{\ell}) -
(\La(y_{\ell})\otimes S(v_\ell)).
$$ 
But this is impossible, since it would give rise to a non-zero element
in $ H^{\ast}(G_\la) - A$ of degree
$d-4\ell +3$, contradicting the minimality of $d$.

\subsection{Computation of $H^{\ast}(G^1_\la)$}\label{ss:1}
\label{algebra-1}

Throughout, we fix $\ell\in\N$ and $\la\in\R$ is such that
$\ell - 1 < \la \leq \ell$. Recall that, for each $0\leq k \leq \ell$,
the strata $U^1_k$ consists of all almost complex structures
$J\in\Jj^1_\la$ for which the class $E_{2k+1}$ has a $J$-holomorphic
representative, and $K^1_k \cong U(2)$ denotes the group of K\"ahler
isometries of a standard complex structure $J_{2k+1}\in U^1_k$.
Moreover, as was shown in \S\ref{description}, $G^1_\la / K^1_k$ is
weakly homotopy equivalent to the strata $U^1_k$.

This whole subsection is devoted to $M^1_\la = (\cp\#\,\ocp, \om^1_\la)$
and the computation of $H^{\ast}(G^1_\la)$, a process which is completely
analogous to what we did in the previous subsection for $H^{\ast}(G^0_\la)$.
To avoid unnecessary repetitions, we will just state the relevant lemmas, 
propositions and corollaries,  leaving all the proofs as exercises 
for the interested reader. In order to simplify notation,
we will also omit the superscript $^1$.

\begin{prop} \label{prop:tensor1}
(i) $H^{\ast}(G_\la)$ is a free algebra.  Moreover
\begin{equation} \label{eq:tensor1}
H^{\ast}(G_\la) \cong H^{\ast}(U_k) \otimes H^{\ast}(K_k)\ ,\
\mbox{for all}\ 0\leq k\leq \ell\ ,
\end{equation}
as graded algebras.

\NI
(ii)  For all $0\leq k\leq \ell$, $H^{\ast}(U_k)\cong H^{\ast}(U_0)$.  
\end{prop}

The computation of $H^{\ast}(G_\la)$
is again based on using the isomorphisms in (i) above and
on the Mayer--Vietoris sequences for the decompositions
$$
U_{0\dots k} =  U_{0\dots \,k-1} \cup \Nn(U_k)
$$
where  $\Nn(U_k)$ is a neighborhood of $U_k$ in $U_{0\dots k}$.  By
Proposition~\ref{prop:link} we may choose $\Nn(U_k)$ to be a disc bundle over
$U_k$, which in particular means that there is a fibration
$$
S^{4k-1}\to \p \Nn(U_k)\to U_k\ .
$$

\begin{lemma}\label{lem:zero1}   For all $ k\leq \ell$
$$
H^i(U_{0\dots k})= 0,\qquad 0< i \leq 4k+2.
$$
\end{lemma}

\begin{proposition}\label{prop:mvk1}
For all $1\le k \le \ell$ there is an exact sequence
$$
(MV_k):
\cdots\to H^i(U_{0\dots k})\to H^i(U_{0\dots
k-1})\stackrel{\tau}{\to} H^{i-4k+1}(U_k)\stackrel{\de}{\to}
H^{i+1}(U_{0\dots k})\to \cdots\ . $$
\end{proposition}

>From $(MV_1)$ and part (ii) of Proposition~\ref{prop:tensor1}, we
get immediately the following

\begin{cor}\label{cor:mv11}
$H^1(U_k) = H^2(U_k) = H^4(U_k) = H^5(U_k) = 0$ and
$H^3(U_k) \cong \Q$, for all $0\leq k\leq \ell$.
\end{cor}

\begin{proposition}\label{prop:deltak1}
(i) In $(MV_k)$, for $k < \ell$, the boundary map
$$
\de_k: \Q = H^3(U_k)\to H^{4k+3}(U_{0\dots k})
$$
is an isomorphism.

\NI
(ii)  If $1\le k < \ell $ the groups $H^{4k+2+j}(U_{0\dots k}) = 0$ for
$j = 2,3$. Moreover
\begin{eqnarray*}
H^{4\ell + 2}(U_{0\dots \ell -1}) & = & \Q,\\
H^{4k+6}(U_{0\dots k}) & = & 0,\quad k < \ell - 1\ .
\end{eqnarray*}
\end{proposition}

Denote by $d$ the smallest integer greater that $5$ such that
$H^d(U_k)\neq 0$. Although independent of $k$ (by 
Proposition~\ref{prop:tensor1}), the value of $d$ depends on the
integer $\ell\in\N$ that was fixed at the beginning of this subsection,
and we make that dependence explicit by writing $d_{\ell}$.

\begin{prop}\label{prop:dl1}  $d_\ell = 4\ell + 2$, for any $\ell\in\N$.
\end{prop}

Since $H^\ast (G_\la)$ is the product of an exterior algebra with a
polynomial algebra, we now know that
$$ \La(a,x,y)\otimes S(w) \subset H^\ast (G_\la)\ ,$$
with $\mbox{deg}\,a=1,\ \mbox{deg}\,x = \mbox{deg}\,y = 3 $ and 
$\mbox{deg}\,w = 4\ell+2$.   Here we define $w$ as in Definition~\ref{def:w}.
Our final task is to show that
these are all the generators of $H^\ast(G_\la)$, i.e., the above inclusion is 
actually an equality of algebras.

For $k=0,1,\ldots,\ell$, let $v_k\in H^{4\ell+2}(U_k)$ denote the
generator corresponding to $w\in H^\ast(G_\la)$ under the isomorphisms
of Proposition~\ref{prop:tensor1}. As in the previous subsection, 
it follows from the proof of the 
previous proposition that the restriction maps
$$ H^i(U_{0\cdots k}) \to H^i(U_{0\cdots k-1}) $$
are isomorphisms when $1\leq k \leq \ell-1$ and $4\ell-2\leq i\leq
4\ell+2$. Hence
$$ i^\ast : H^{4\ell+2}(U_{0\cdots k}) \to H^{4\ell+2}(U_0) $$
is an isomorphism for $0\leq k\leq \ell-1$, and we also denote by
$v_0$ the generator of each $H^{4\ell+2}(U_{0\cdots k})$ such that
$i^\ast(v_0)=v_0$.

\begin{lemma} \label{lem:ik1} Let $i_k$ denote the inclusion
$U_k \to U_{0\cdots k}$. Then  $i_k^\ast (v_0) = v_k$, for $0\leq k\leq
\ell-1$. \end{lemma}

The collapse of the spectral sequence of the fibration
$$ S^{4k-1}\to\p\Nn(U_k)\stackrel{\pi}{\to} U_k\ ,$$
together with the fact that the fiber is an odd-dimensional sphere,
implies that we have ring isomorphisms
$$H^{\ast}(\p\Nn(U_k))\cong H^{\ast}(S^{4k-1})
\otimes H^{\ast}(U_k),\ 1\leq k\leq\ell\ .$$
Hence we know that
$$H^{\ast}(\p\Nn(U_k))\supseteq \La(\pi^\ast(y_k),
\widehat{z}_{4k-1})\otimes S(\pi^\ast(v_k)),\ 1\leq k\leq\ell,$$
where the elements $y_k  \in H^{3}(U_k)$ and 
$\widehat{z}_{4k-1}\in H^{4k-1}(\p\Nn(U_k))$
are  as before.  Thus the $y_k$ are 
are defined as in
Definition~\ref{def:w}, and  $\widehat{z}_{4k-1}$ is the restriction of $z_{4k-1}\in
H^{4k-1}(U_{0\cdots k-1})$, where $z_{4k-1}$ is determined from the
isomorphism $\tau: H^{4k-1}(U_{0\cdots k-1})\to H^{0}(U_k)$ in $(MV_k)$ by
the condition $\tau(z_{4k-1}) = 1$.

\begin{lemma} \label{lem:rest1} Let $j_k$ denote the inclusion 
$\p\Nn(U_k) \to U_{0\cdots k-1}$ and $\pi$
the projection $\p\Nn(U_k)\to U_k$. Then
$j_k^{\ast}:H^{\ast}(U_{0\cdots k-1})\to H^{\ast}
(\p\Nn(U_k))$ satisfies:
\begin{itemize}
\item[{\rm (i)}] $j_k^{\ast}(v_0) = \pi^{\ast}(v_k)$, for $0\leq k\leq\ell-1$;
\item[{\rm (ii)}] $j_\ell^{\ast}(v_0) = \widehat{z}_{4\ell - 1}\cdot
\pi^{\ast}(y_\ell) + c \pi^{\ast}(v_\ell)$, with $c\neq 0$.
\end{itemize}
\end{lemma}

\begin{lemma} \label{lem:free1}
$z_{4k-1}$ and $v_0$ generate free subalgebras of 
$H^{\ast}(U_{0\cdots k-1})$, for $1\leq k\leq\ell$.
\end{lemma}

\begin{proposition} \label{prop:maps1}
In the exact sequence $(MV_{k})$
$$\cdots\to H^i(U_{0\cdots k})\stackrel{i^{\ast}}{\to} H^i(U_{0\cdots
k-1})\stackrel{\tau}{\to} H^{i-4k+1}(U_k)\stackrel{\de}{\to}
H^{i+1}(U_{0\cdots k})\to\cdots$$
(for $1\leq k\leq\ell-1$), the maps $i^{\ast}$, $\tau$ and $\de$ 
satisfy:
$$
\begin{array}{lll}
i^{\ast}(v_{0}^n) = v_{0}^n & \tau(v_{0}^n) = 0 & \de(v_{k}^n) = 0 \\
i^{\ast}(z_{4k+3}\cdot v_{0}^n) = 0 & \tau(z_{4k-1}\cdot v_{0}^n)
= v_{k}^n & \de(y_{k}\cdot v_{k}^n) = c\cdot z_{4k+3}\cdot v_{0}^n,\ 
c\neq 0\ .
\end{array}
$$
\end{proposition}

The final argument to show that $$H^{\ast}(G_{\la}) = \La (a,x,y) \otimes S(w)$$
is exactly the same as in the previous subsection. More precisely, 
if $$H^{\ast}(G_{\la}) \neq A = \La (a,x,y) \otimes S(w),$$ then
we can choose a non-zero element $u$ of minimal degree $d$ in
$
H^{\ast}(G_{\la}) - A,
$
and show that this new element $u$ would give rise to 
a new non-zero element in $H^{\ast}(G_\la) - A
$ of degree $d-4\ell +1$, which by
minimality of $d$ is a contradiction.

This completes the proof of Theorem~\ref{thm:main1}.

\section{The limit at infinity }\label{sec:lim}

In this section we consider the limit structure of the groups $G_\la^i$ as
$\la \to \infty$.  Our aim is to show that this limit exists and can be
identified (up to homotopy) with the group $\Dd$ of fiberwise
diffeomorphisms considered in \S\ref{ss:Kk}.
 To begin with we will consider only the case $i = 0$
and, as usual, will omit superscripts.  However, as we will see in
Remark~\ref{rmk:any}, our arguments apply to any ruled surface.

In order to consider the  limit $\lim_{\la\to
\infty}G_\la$,  we need  to construct maps $G_\la\to
G_{\la+\ka}$ for $\ka>0$.  However, there is no very direct way to do this,
and the best approach is rather to include $G_\la$ into a space
$G_{[\la,\la+\ka]}$ that deformation retracts onto $G_{\la+\ka}$.
Therefore we  make the following definitions.

\begin{defn}\label{def:inf}\rm
Given a family $H_\la, \la \in [0,\infty)$  of subgroups of $\Diff(S^2\times
S^2)$ that varies smoothly with $\la$, define
\begin{eqnarray*}
 H_\infty = \bigcup_{\la> 0} \{(\phi,\la): \phi\in H_\la\}\subset
\Diff(S^2\times S^2)\times \R^+,\\
H_{[\la,\la+\ka]} = \bigcup_{\la\le \mu\le \la+\ka} \{(\phi,\mu): \phi\in
H_\mu\}\subset \Diff(S^2\times S^2)\times \R^+. 
 \end{eqnarray*}
\end{defn}

Our first aim is to prove the following result, which shows that
$$
G_\infty = \lim_{\la\to \infty} G_{[0,\la]}
$$
is the homotopy limit of the $G_\la$ as $\la\to \infty$.

\begin{prop}\label{prop:equiv}  The inclusion $G_{\la+\ka}\to G_{[\la,\la+\ka]}$ is a
weak homotopy equivalence.
\end{prop}

To prove this proposition we must consider various enlargments of the groups
$G_\la$.

\begin{defn}\rm  We define $\TG_\la$ to be the set of all pairs $(\phi,
\tau_t)$, where $\phi\in \Diff^+(S^2\times S^2)$ and $\tau_t, 0\le t\le 1$, is
a smooth family of cohomologous symplectic forms such that $\tau_t =
\om_\la$ for $t$ near $0$ and $\tau_t = \phi^*(\om_\la)$ for $t$ near $1$.
Define an operation $\TG_\la\times \TG_\la \to \TG_\la$ by
$$
(\phi, \tau_t)\circ(\phi', \tau_t') = (\phi\circ\phi', (\tau*\tau')_t)
$$
where 
$$
 (\tau*\tau')_t   = \left\{\begin{array}{lll} \tau_{2t}' & \mbox{if} & t\in [0,1/2]\\
 (\phi')^*\tau_{2t-1} & \mbox{if} & t\in [1/2,1]
\end{array}\right.
$$
With this operation, $\TG_\la$ is a homotopy associative $H$-space with
homotopy identity and homotopy inverse. 
Moreover, $G_\la$ can be identified with the subgroup of  $\TG_\la$ 
on which the path $\tau_t$ is constant.
\end{defn}

The next lemma follows immediately from the usual Moser argument.

\begin{lemma} $\TG_\la$ deformation retracts onto its subgroup $G_\la$.
\end{lemma}

We would like to  define a family of maps $G_\la \to \TG_{\la + \ka}$.
However, as we shall see, we have to do something a little more elaborate.
The next lemma will be useful.

\begin{lemma}\label{le:nd}  Let $\pi: X\to \Si$ be a smooth fibration of an
oriented $4$-manifold $X$ over an oriented $2$-dimensional base manifold, and
suppose that $\tau$ is a symplectic form on $X$ whose restriction to each fiber
is nondegenerate. Then, if $\be$ is any nonnegative form on $\Si$, the form $\tau
+\pi^*(\be)$ is symplectic.
\end{lemma}
\proof{} For each $x\in X$ let $H_x\subset T_xX$ be the $2$-dimensional
$\tau$-orthogonal complement to the tangent space $\ker d\pi(x)$ to the fiber
at $x$. An easy calculation shows that a form $\rho$ must be nondegenerate if it
 restricts to $\tau$ on the fibers and is a positive multiple of $\tau$ on
each $H_x$.\QED

Let $\Ss_\la$ denote the space of symplectic forms on $S^2\times S^2$ that
are isotopic to $\om_\la$, 
 $\Fol_\la$ be the space of smooth foliations of
 $S^2\times S^2$ whose leaves are $\om^\la$-symplectic spheres in class
$F$,
and  $\Proj_\la$ be the space of smooth surjections $p:S^2\times S^2\to S^2$
whose fibers form an element of $\Fol_\la$ and that are compatible with the
obvious orientations.   Then
there is a fibration
$$
\Diff^+(S^2) \to \Proj_\la \to \Fol_\la.
$$
Moreover, $\Diff^+(S^2)$ acts on the left on $\Proj_\la$ via 
$ \psi\cdot p = \psi\circ p$, and the induced map
$$
\SO(3) \!\setminus\! \Proj_\la \to \Fol_\la
$$
is a homotopy equivalence.  The preceding lemma implies that there is a 
family of maps
$$
s_\ka: \Proj_\la \to \Ss_{\la+\ka}
$$
given by
$$
s_\ka (p) = \om_\la + \ka p^*(\si),
$$
where $\si$ is the usual $\SO(3)$-invariant area form on $S^2$ of total area
$1$.  Since $\si$ is $\SO(3)$-invariant, we in fact have a map
$$
s_\ka: \SO(3)\!\setminus\! \Proj_\la \to \Ss_{\la+\ka}.
$$

\begin{lemma} The map
$$
\Psi: \Jj_\la \to \Fol_\la$$
that takes $J$ to the foliation by $J$-holomorphic spheres in class $F$
is a homotopy equivalence.
\end{lemma}
\proof{}  This holds by standard arguments because the inverse image of
each point in $\Fol_\la$ is contractible.\QED

Each of the spaces $\Ss_\la, \Fol_\la, \Proj_\la$ has an obvious base point.
If $Y$ is any based space we write $\Pp_*(Y)$ for the space of smooth paths
in $Y$ starting at the base point.

The choice of a deformation retraction
of $\Jj_\la$ onto the split complex structure $J_0$ defines a 
map $\io_r: \Jj_\la \to \Pp_*(\Jj_\la)$ and hence a map
$$
G_\la \to \Pp_*(\Fol_{\la}) : \phi\mapsto  \phi^*(J_0)
\mapsto \io_r(\phi^*(J_0)\mapsto \Psi_*(\io_r(\phi^*(J_0))).
$$
Thus there is a diagram
$$
\begin{array}{ccccc}
& & \Pp_*(\SO(3)\!\setminus\! \Fib_\la) & \stackrel{s_\ka}\to &
\Pp_*(\Ss_{\la+\ka})\\
&f\nearrow & \downarrow\simeq & &\\
G_\la & \longrightarrow & \Pp_*(\Fol_{\la}) & & 
\end{array}
$$
in which the lifting $f: G_\la \to \Pp_*(\SO(3)\!\setminus\! \Fib_\la)$ is really
defined only over compact subsets of $G_\la$.
Thus for each compact subset $C$ of $G_\la$ we have a map
$$
f_C: C\to\TG_{\la+ \ka}:\quad \phi\mapsto (\phi, s_\ka(f(\phi)).
$$
Moreover, this set of maps $f_C$ is compatible in the sense that if
$C_1\subset C_2$ then the restriction of $f_{C_2}$ to $C_1$ is homotopic to 
$f_{C_1}$.
\MS\MS

\NI
{\bf Proof of Proposition~\ref{prop:equiv}} 
 This follows immediately from the preceding discussion. \QED

\MS

Our next aim is to show that $G_\infty$ is homotopy equivalent to 
the group $\Dd$ of orientation preserving fiberwise diffeomorphisms,
which up to homotopy is the same as
$$
\Dd = \SO(3)\times \Map(S^2, \SO(3))\ .
$$.
Since the elements $\phi\in \Dd$ preserve the fibers of the standard
fibration $\pi:S^2\times S^2\to S^2$, Lemma~\ref{le:nd} implies that the
forms $\phi^*(\om_\la)  + \ka\pi^*(\si)$ are symplectic for all $\ka$. 
Moreover, for sufficiently large $\ka$ the elements of the linear isotopy
$$
t(\phi^*(\om_\la)  + \ka\pi^*(\si)) + (1-t)\om_{\la+\ka}, \quad 0\le t\le 1,
$$
are also all symplectic.  Thus, if
  $\Dd_K$ is the subset of $\Dd$ consisting of elements
whose first derivative is uniformly bounded by $K$, there is, for some 
constant $\ka(K)$ a
natural map $$
\Dd_K \to \TG_{\la + \ka},\quad \ka \ge \ka(K).
$$
Hence
there is a map from
the homotopy limit 
$$
\Dd' = \lim_{\rightarrow}  \Dd_K
$$
to $\TG_\infty$.   The obvious map $\Dd'\to \Dd$ is a weak homotopy
equivalence, and hence a homotopy equivalence since both spaces have the
homotopy types of CW complexes.  Thus we may 
 define a map $\Dd\to G_\infty$ via the composite
$$
\Dd\to \Dd' \to \TG_\infty \to G_\infty.
$$

To get a map the other way round, we must enlarge the group $\Dd$.
For each $\la$, let 
$$
\Dd_{\la}^{\Proj} = \{\phi\in \Diff(S^2\times S^2): \pi \circ\phi \in \Proj_\la\}. 
$$
It is not hard to check that there is a fibration 
$$
\Dd\to \Dd_{\la}^{\Proj}  \to {\Proj_\la}.
$$
(One can construct
local trivializations of $\Dd_{\la}^{\Proj}  \to {\Proj_\la}$ near a 
projection map $p$ using
 a section of $p$.)  Hence, since $\Proj_\la$ is contractible, $\Dd\simeq
\Dd_{\la}^{\Proj}$ for all $\la$.  Thus, defining 
$\Dd_{\infty}^{\Proj}$ as in Definition~\ref{def:inf}, we have
$$
\Dd
\simeq \Dd_{\infty}^{\Proj}
$$
\MS

Since
$G_\infty\subset \Dd_{\infty}^{\Proj}$, we can define a map $G_\infty\to \Dd$
as the composite
$$
G_\infty\to \Dd_{\infty}^{\Proj} \stackrel{\simeq}\longrightarrow \Dd.
$$
This is the homotopy inverse to the map $\Dd\to G_\infty$ defined before
because, by Lemma~\ref{le:nd}, every compact subset of $\Dd_{\la}^{\Proj}$ is
contained in $\TG_{\la + \ka}$ for suitably large $\ka$.

\MS

\begin{remark}\label{rmk:any}\rm  The analog of this result is true for any
ruled surface $X$, i.e. the limit, when the ratio  of the size of the base
to that of the fiber goes to infinity, of the group of symplectomorphisms of
$X$ is homotopy equivalent to the group of fiberwise diffeomorphisms
$\Dd(X)$ of $X$ that act trivially on $H_2(X)$.  
When $X$ is the product
$\Si\times S^2$, this group $\Dd(X)$ has the homotopy type of the product
$\Diff^+(\Si)\times \Map(\Si,\SO(3))$.
 In the case of the nontrivial bundle,
 let $\Ga(\Si;\SO(3))$
denote the space of sections of the bundle over $\Si$ with fiber $\SO(3)$
that is associated to $X$.  It is not hard to check that every area-preserving
diffeomorphism $\phi_\Si$ of $\Si$ does lift to $X$: if $\phi_\Si$
is the time $1$ map of the Hamiltonian $H_t$, a lift is given by the time
$1$ map of $H_t\circ \pi$.   
When $X\to S^2$ is the nontrivial $S^2$-bundle over $S^2$ we worked out the
homotopy type of the relevant group $\Dd(X)$ in \S\ref{ss:Kk}: namely
$$
\Dd(X)\simeq\SU(2)\times_{\Z/2\Z} \Ga(S^2;\SO(3)).
$$
 When $genus(\Si) > 0$ the fibration
$X\to \Si$ is pulled back from $\cp\#\,\bcp\to S^2$, and it is not
hard to see that elements in all components of $\Diff^+(\Si)$ lift.  Moreover, 
the action of $\Diff^+(\Si)$ on $\Ga(\Si;\SO(3))$ is trivial.  
Hence, in this case
$$
\Dd(X)\simeq \Diff^+(\Si)\times \Ga(\Si;\SO(3)).
$$

One can also consider the other limit, when the relative size of the base goes to zero.
The case $X = S^2\times S^2$ is exceptional, since there is an extra symmetry of
$X$ that interchanges the base and fiber.  For all other $X$ the fiber class $F$ is
distinguished by the fact that  every tame $J$ defines a fibering of $X$ by
$J$-holomorphic spheres in class $F$.    When $X = \cp\#\,\ocp$ the
limit that we are now considering is given by $\la \to -1$, and the groups
$G_\la$ stabilize: by Gromov's work they all equal $U(2)$ as soon as $\la
\le 0$.  It would be interesting to know if this result had analogs  when $genus(\Si)
> 0$. \end{remark} 

\MS

\section{Relation with Kronheimer's work}
\label{ss:kron}

In~\cite{K}, Kronheimer described some nontrivial families $\Ff_X$ of 
cohomologous symplectic
structures on certain $4$-manifolds $X$.  In particular he showed how to
construct a family $\Ff_\ell$ of  symplectic forms $\om^\nu$ on $
S^2\times S^2$ in class $[\om_\la]$ that are parametrized by
$\nu\in S^{4\ell+1}$.
In many cases,  he  shows
by an argument using Seiberg--Witten invariants that 
$\Ff_X$ represents a nontrivial
homotopy class in the appropriate space $\Ss_X$ of 
cohomologous symplectic forms on $X$.  Since his argument only works 
 when $b^+(X)$ is sufficiently large, it does not apply directly to $S^2\times S^2$
or to $\cp\#\,\bcp$.
 In this section we exploit the close connection between $\Ff_\ell$ and  $w_\ell$
to show that $\Ff_\ell$ is nontrivial also in these cases.

As was pointed out at the beginning of the introduction, both 
$M^{0} = S^2\times S^2$ and $M^{1} = \cp\#\, \ocp$ have a unique 
symplectic structure in each permissible cohomology class. In 
particular, this means that the spaces 
$$
\Ss_{\la}^{i} \equiv \{ \mbox{all symplectic forms on}\ M^{i}\ 
\mbox{representing the cohomology class}\ [\om_{\la}^{i}]\} 
$$
are homogeneous for the groups $\Diff(M^{i})$ of all diffeomorphisms of 
$M^{i}$:
$$
 \Ss_{\la}^{i} \cong \Diff(M^{i})/G_{\la}^{i},\ i=0,1.$$
The corresponding fibration
$$
G_{\la}^{i} \longrightarrow  \Diff(M^{i}) 
\longrightarrow \Ss_{\la}^{i}
$$
has an associated long exact sequence
$$
\cdots \to \pi_{k}(\Diff(M^{i}) ) \to \pi_{k}(\Ss_{\la}^{i}) 
\stackrel{\p}{\to} \pi_{k-1}(G_{\la}^{i})
\to \pi_{k-1}(\Diff(M^{i}) ) \to \cdots \ .
$$
Since the rational cohomology of any $H$-space is freely generated by its
rational homotopy,  Theorems~\ref{thm:main0} and~\ref{thm:main1}
imply that $\pi_*(G_\la^i)\otimes\Q$ has $4$ generators, in degrees $1$, $3$,
$3$ and $4\ell + 2i$.  The first three of these map to nonzero elements in
$\Diff(M^i)$,
but the fourth does not since it already vanishes in $G_{\la + 1}^i$.
Hence it must lie in the image of the boundary map $\p$ of the above
long exact sequence.   Thus we have:

\begin{prop} \label{prop:sympforms}
For $i=0,\,1$, and $0 \le \ell -1<\la\leq \ell,\ \ell\in\N$, there exists a 
non-zero element  $\Ff_\la \in\pi_{4\ell + 1 + 2i}(\Ss_{\la}^i )$ 
whose image in $\pi_{4\ell + 2i}(G_{\la}^i)$ 
maps to a non-trivial generator of $H_{4\ell + 2i}(G_{\la}^i , \Q)$.
\end{prop}

\begin{remark} \rm When $\ell = 0$ (i.e. $\la =0$ for $i=0$, or 
$-1<\la\leq 0$ for $i=1$) the boundary map $\p:\pi_{1 + 2i}(\Ss_{\la}^i )
\to \pi_{2i}(G_{\la}^i)$ is zero. 
\end{remark}

 Kronheimer gave an explicit construction for this element $\Ff_\ell$
at least in the case $\la = \ell$, which shows clearly why it vanishes when
$\la$ increases past $\ell$.  We will first describe this construction and 
then relate it to the picture given  in Definition~\ref{def:w} of
$w_\ell$.  As usual, we restrict to the case $i = 0$ and omit superscripts.  The
case $i = 1$ is almost identical.
\MS

The construction goes as follows.  Let $V$ be an open neighborhood 
of $\{0\}\in \C^{2\ell + 1}$.  For
$\ell\ge 1$ Kronheimer  constructs in~\cite{K}~\S3 a smooth family 
$$
(X^v, J^v,
\tau_\ell^v),\qquad v\in V - \{0\},
$$
 of  K\"ahler manifolds diffeomorphic to
$X = S^2\times S^2$ with K\"ahler  form $\tau_\ell^v$ in class $[\om_\ell]$. 
This family can be  completed by a complex fiber $(X^0,J^0)$ over $0$ in 
such a way that the following conditions are satisfied.
\begin{itemize}
\item[(i)]  The whole family ${\bf X} = \cup_{v\in V}X^v$, is smoothly 
trivial, i.e. is fiberwise diffeomorphic to $X\times V$. Moreover there is a
complex structure on ${\bf X}$ that restricts to $J^v$ on each fiber $X^v$, and
a closed $(1,1)$ form
$\tau_\ell$ on ${\bf X}$ that restricts to
$\tau_\ell^v$ when $ v\ne 0$. 
\item[(ii)] The restriction of   $\tau_\ell$ to $X^0$
degenerates along a holomorphic curve $\Delta$ in class
$A -  (\ell+1)F$. 
\item[(iii)] For any $\la > \ell$ there is a  K\"ahler 
form $\tau_\la$ on ${\bf X}$ whose restriction $\tau_\la^v$ to each $X^v$
lies in the class $ [\om_\la]$. Moreover the family $\tau_\mu, \mu\ge\ell,$
varies smoothly with $\mu$.
\end{itemize}

Condition (ii) implies that   $(X^0, J^0)$ is simply $ (S^2\times S^2, J_{2\ell+
2})$.   The complex variety $Y^0$ that is obtained from $X^0$ by collapsing
the curve $\Delta$ has  a quotient singularity of the form $\C^2/C_{\ell+1}$
where $C_p$ is the cyclic group of order $p$ acting diagonally.
If ${\bf Y}$ is obtained from  ${\bf
X}$ by collapsing $\Delta$, ${\bf Y}$ is smooth and the form $\tau_\ell$
descends there.  In fact, as is explained in~\cite{K}~\S5, ${\bf Y}$ is 
the Artin component of the universal unfolding of this 
given quotient singularity, and ${\bf X}$ is 
its simultaneous resolution.
\MS

Several points in this construction are relevant to us.
First observe that we can consider the 
forms $\tau_\la^v$ to lie on the fixed manifold $X = S^2\times S^2$ by 
using the smooth trivialization given in (i).  When $\la = 
\ell$ we therefore  get a family $\tau_\ell^\nu, \nu\in V - \{0\},$ of 
K\"ahler forms on $X$ and hence an element $\Ff_\ell'\in \pi_{4\ell +
1}(\Ss_\ell)$.   We will show:

\begin{prop}\label{prop:bound} $\p(\Ff_\ell') \ne 0$ where   $\p:
\pi_{*+1}(\Ss_\ell)\to \pi_*(G_\ell)$ is the boundary map considered above.
\end{prop}

Thus we can take $\Ff_\ell = \Ff_\ell'$ in Proposition~\ref{prop:sympforms}.
This picture explains why $\Ff_\la$ vanishes when $\la$ increases beyond
$\ell$.  For, as mentioned in (iii) above,  the family
$\tau_\ell^v, v\in V - \{0\},$ smoothly deforms to a family
$\tau_\la^v,$
 $\la > \ell,$ that extends to $v = 0$.  Using this, it is easy to see that if
$G_\ell\to G_\la$ is constructed as in \S\ref{sec:lim}, $\p(\Ff_\ell)$
is in the kernel of the
induced map on $\pi_{4\ell}$.\MS

As a first step in the proof of this proposition, observe that when $\la > \ell$ we can
trivialize the family $
(X^v,\tau_\la^v), v\in V,
$
 as a family of {\it symplectic} manifolds,
identifying them all with $(S^2\times S^2, \om_\la)$.   This gives us a family
$$
\Vv_\la = \{J^\nu: \nu\in V\}
$$
 of elements of $\Jj_\la$.  Condition
(ii) implies that $J^0$  is an element in the stratum $U_{\ell+1}$
that we can identify with $J_{2\ell + 2}$. 
Moreover, as pointed out at the end of~\cite{K}~\S4, this
family $\Vv_\la$ is transverse to  $U_{\ell+1}$.  
In fact, since $U_{\ell+1}$ has (real) codimension precisely $4\ell + 2$ in $\Jj_\la$,
$\Vv_\la$ is a normal slice to $U_{\ell+1}.$
Therefore, we can investigate the structure of the links of one stratum in
another by looking at the intersections of these strata with $\Vv_\la$.
In particular, none of these intersections is empty and the
boundary of the intersection $U_\ell\cap \Vv_\la$ 
can be identified with the link
$L_{\ell+1}$.  

\MS

Now, let us look at the corresponding picture when $\la = \ell$.  
In this case, we restrict $v$ to a small $(4\ell + 1)$-sphere  $S_V$ that
encloses $\{0\}$ in $V$ and pick a basepoint $*\in S_V$ at which
$ J^{*}$ is diffeomorphic to the standard split structure. 
Let $D$ be a neighborhood of $*$ in $V_S$ that is
diffeomorphic to an open disc and is small enough that the $J^v$ for 
$v\in {\ov{D}}$  are all split.  Then trivialize  the symplectic manifolds
$(X^v,\tau_\ell^v)$ for $v$ in the disk $V_S - D$, identifying them all with
$(S^2\times S^2, \om_\ell)$. This gives us a smooth family
$$
\Ww_\ell = \{J^\nu: \nu \in V_S - D\}\subset \Jj_\ell.
$$
Moreover, our choice of $D$ implies that $\p\Ww_\ell$ lies entirely in
$U_0$.  Thus we get an element $\be_\ell$ of $\pi_{4\ell}(U_0)$.  

\begin{lemma}  This element $\be_\ell\in \pi_{4\ell}(U_0)$ lifts to
$\p(\Ff_\ell')\in \pi_{4\ell}(G_\ell)$.
\end{lemma}
\proof{}  To see that $\be_\ell$ lifts to $G_\ell$, we just need to perform the
previous construction with a little more care.
By construction, we can identify all the K\"ahler manifolds $(X^v, J^v,
\tau_\ell^v),$ $v\in {\ov{D}},$ with the standard model $(S^2\times S^2, J_0,
\om_\ell)$.    Then the manifolds $(X^v,\tau_\ell^v), v \in \p D,$ have two
identifications with $S^2\times S^2$ and the difference between these two
gives the desired lift.  This lift represents $\p(\Ff_\ell')$ by the definition of the
boundary map $\p$.\QED

The next claim is that the sets
$$
\Ww_\ell\subset \Jj_\ell,\quad\mbox{and}\quad  \Vv_\la \subset \Jj_\la,\la >
\ell, $$
fit together nicely as $\la$ varies.  More precisely,
define
$$
\Ww_\la = \{J^\nu: \nu\in S_V - D\},
$$
and, for some $\eps > 0$, consider the space
$$
\Ww = \bigcup_{\mu\in [\ell, \ell + \eps)} \,\Ww_\mu.
$$
Each space $\Ww_\la$ inherits a stratification from $\Jj_\la$, and it is not
hard to see that by careful use of the Moser method one can construct all the
symplectic trivializations  to vary smoothly with $\mu$.   Thus we have:

\begin{lemma}  The stratifications of each $\Ww_\mu$ fit together to give
a stratification of $\Ww$ by the sets $U_0\cap\Ww,\dots, U_\ell\cap \Ww$. 
In particular,  the
intersections $U_k\cap \Ww_\mu$ for $\mu = \ell$ and $\mu > \ell$ are
diffeomorphic. 
\end{lemma}

\begin{cor}  The intersection $U_\ell\cap \Ww_\ell$ is a 
connected $3$-dimensional set that represents a  generator of $H_3(U_\ell)$.
\end{cor}
\proof{} By the previous lemma it suffices to prove this for $\mu > \ell$.  But we 
remarked before that the boundary of the intersection $U_\ell\cap \Vv_\la$ 
can be identified with the link
$L_{\ell+1}$.  Now use Corollary~\ref{cor:hlink}.\QED

\begin{remark}\rm
One can construct $\Ww_\ell$ so that it is invariant under the action of 
the $\SU(2)$ factor in the
automorphism group $K_{\ell+1}$ of $J_{2\ell+2}$.
Then the intersection $U_\ell\cap \Ww_\ell$ is just an $\SU(2)$-orbit 
and so, by Lemma~\ref{le:htpyK0}, represents the class $\xi + 
(\ell+1)^2 \eta$.
\end{remark}

The proof of Proposition~\ref{prop:bound}   is completed by the
following result.

 \begin{prop} The image $v_0$ of $w_\ell$ in
$H^{4\ell}(U_0)$ does not vanish on $\be_\ell$.  
\end{prop}
\proof{}  $\Ww_\ell$ is a  $(4\ell + 1)$-dimensional disk in $\Jj_\ell$
that intersects all the strata $U_k$ transversally.  
We have just seen that it intersects $U_\ell$ in a cycle $C$ that generates
$H_3(U_\ell)$.  Let $Z\subset \Ww_\ell$ be the boundary of a tubular 
neighborhood $\Nn(C)$ of $C$.  Then, $Z$ has dimension $4\ell$ and lies in
$
(U_{0\dots\ell-1})\cap \Ww_\ell.
$
In the notation of Remark~\ref{rmk:w},  $Z = C*D$, where in fact
$C = \xi+(\ell+1)^2 \eta$ (though this will not be needed.)
Thus ${\widehat{y}}_\ell(Z) \ne 0$.   Further, since $C = \Ww_\ell\cap U_\ell$,
$$
\Ww_\ell - \Nn(C) \subset U_{0\dots \ell -1}.
$$
Since $C$ is connected, $H^{4\ell -4}(\Ww_\ell - \Nn(C)) = \Q$.  Hence the cycle $Z -
\be_\ell$ is null homologous in $\Ww_\ell - \Nn(C)$ and hence in $U_{0\dots \ell
-1}$. Thus
$$
 v_0 (\be_\ell)= {\widehat{y}}_\ell (\be_\ell) = {\widehat{y}}_\ell(Z)\ne 0,
$$
as claimed.\QED

\section{Whitehead products and $H^*(BG_\la)$}
\label{ss:whit}

Here we exhibit 
the $w_\ell$ as dual to certain higher  products
in the Lie algebra $\pi_*(G_\la)\otimes\Q$.  Thus they are 
 desuspensions of higher Whitehead products in $\pi_*(BG_\la)$, and hence,
via the work of Andrews--Arkowitz~\cite{AA},
give rise to  relations in
 the rational cohomology of the classifying spaces
$H^*(BG_\la)$.   In fact, our knowledge of $H^*(G_\la)$ allows us to calculate
the differential graded Lie algebra $\pi_{*-1}(BG_\la)\otimes \Q$ and hence 
to determine the rational homotopy type of $BG_\la$.\footnote
{The calculations in~\S\ref{sec:algebra} show that $G_\la$, and hence $BG_\la$, 
has finite rational homotopy type.  Moreover $\pi_1(BG_\la)\otimes \Q  = 0$.
Thus the work of~\cite{AA} applies.}
  In particular, since $H^*(BG_\la)$ is not a free ring,
$BG_\la$ does not have the rational homotopy type of an $H$-space.
\MS

 First recall (see for example Ch X of~\cite{Wh}) that
for any group $G$ the Samelson product $[\al,\be]\in \pi_{p+q}(G)$ of
elements  $\al\in\pi_p(G)$ and $\be\in \pi_q(G)$
is represented by the map
$$
S^{p+q} = S^p\times S^q/ S^p\vee S^q \quad\longrightarrow X:
(u,v)\mapsto \al(u)\be(v)\al(u)^{-1}\be(v)^{-1}.
$$
(For simplicity of notation, we will use the same letter for a map
$\tau:C\to X$ as for the element of homotopy or homology that it
represents.)
The Samelson product in $\pi_*(G)$ is related to the Pontrjagin product
$*$ in $ H_*(G, \Z)$ by the formula
$$
[\al,\be] = \al*\be - (-1)^{pq} \be *\al\quad\in H_*(G).
$$
Here we have suppressed mention of the Hurewicz map
$\pi_*(G)\to H_*(G)$, and have written $\al*\be$ for the cycle
$$
S^p\times S^q\to G:\quad (u,v)\mapsto \al(u)\cdot \be(v).
$$
Thus the domain $\Dom(\al*\be)$ of $\al*\be$ is $S^p\times S^q$, while 
$\Dom(\be*\al)$ is $S^q\times S^p$.   Moreover,
if we think of $\al,\be\in \pi_*(G)$ as the 
desuspensions of elements $E(\al),
E(\be)$ in $\pi_*(BG)$, the Samelson product $[\al,\be]$ is (up to sign) the
desuspension of the Whitehead product $[E(\al), E(\be)]\in \pi_{p+q+1}(BG)$:
see~\cite{Wh} Thm~(X.7.10).  Hence the Samelson products in $G$ vanish 
if $BG$ itself is  an $H$ space, since Whitehead products vanish in $H$-spaces.

We will say that  $\al,\be$ are {\it commuting representatives} for their
respective classes if 
$$
\al(u)\be(v) = \be(v)\al(u) \quad\mbox{for all }\;\; u,v\in
\Dom(\al)\times \Dom (\be).
$$
  Thus, in this case, $[\al,\be]$ is trivially $0$.

Let us first consider the case $i = 0$ and $0 < \la\le 1$.  As usual we
suppress the superscript $i$ and write $w$ instead of $w_\ell$.
We will also write $G$ instead of $G_\la$ and will 
work with rational coefficients.
Further, let  $\al = \al_1:S^1\to G$ and $\xi= \xi_0:S^3\to G$ be as
in~\S\ref{ss:Kk}.  
  Note that $[\xi,\eta] = 0$ because $\xi,\eta$ have commuting
representatives in $K_0$, while $[\xi,\xi]$ and $[\eta, \eta]$ vanish
  because 
$\pi_6(S^3)$ is finite.  
\MS

\begin{lemma}\label{le:sa1} For $0 < \la\le 1$, $[\al,\xi+\eta] = 0$ while
$w([\al,\xi]) \ne 0$. \end{lemma}
\proof{}
By Lemma~\ref{le:htpyK0} the classes $\al$ and $\xi+\eta$ have commuting
representatives in $K_1\subset G =G_\la$, which proves the first statement.

To prove the second, note first that because
 $\xi$ has image in $K_0$, it maps to a trivial
element in $U_0 = G/K_0$.  Hence the image of $[\al,\xi] = \al*\xi + \xi *\al$
in $H_4(U_0)$ is simply the Pontrjagin product $\xi*\al$.  Now $G$  maps to
$U_0$  via $\phi\mapsto \phi_*(J)$ where we can choose any $J\in U_0$.  If
we choose $J\in \Nn(U_1)$ very close to $J_1\in U_1$ then the loop
$\th\mapsto \al(\th)_*(J)$ circles the stratum $U_1$ while the sphere
$u\mapsto \xi(u)_*(J)$ is a copy in $\Nn(U_1) - U_1$ of a generating cycle for
$H_3(U_1)$.  Thus the map 
$$
H_*(\p\Nn(U_1))\to H_{*-1}(U_1)
$$
induced by integration over the fiber takes the Pontrjagin product $\xi*\al$
to a generator of $H_3(U_1)$.  Since this map is dual to the map $\tau$
occuring in $(MV_1)$ it follows from Definition ~\ref{def:w}
that the element $v_0\in H^4(U_0)$ that corresponds to $w$ takes a nonzero
value on $\xi*\al$.\QED

\begin{cor}\label{cor:BG1}When $0< \la\le 1$, $H^*(BG_\la) =
\Q[A,X,Y]/\{A(X-Y) = 0\}$, where $A$ has dimension $2$ and $X,Y$  have
dimension $ 4$. \end{cor}
\proof{}  The structure theorem for the rational cohomology of an
$H$-space $G$ states that $H^*(G)$ is freely generated as a graded algebra by
the rational homotopy groups $\pi_*(G)\otimes \Q$.  Thus $\pi_*(G)\otimes
\Q = 0$ for $*\ne 1,3,4$, and is generated by the elements $\al,\xi,\eta,
[\al,\xi]$.  
 Since $
[\al,\xi+\eta] = 0,
$
 the only nontrivial  Whitehead product in $BG$ (including higher
order ones) is the surjective map
$$
[\cdot, \cdot]: \pi_2(BG)\otimes \pi_4(BG)\otimes \Q \to \pi_5(BG)\otimes \Q = \Q,
$$
with kernel $E(\al)\otimes \Q(E(\xi+\eta))$. By~\cite{AA} Theorem~5.4,
the dual of this map corresponds to a relation in $H^*(BG)$.

More precisely, when one constructs the minimal model $\Mm$ of $H^*(BG)$,
there is a generator $A\in \Mm^2$ dual to $E(\al)$, generators $X$ and $Y$ in
$\Mm^4$ dual to $E(\xi)$ and $ E(\eta)$, and a generator $W$ in $\Mm^5$
corresponding to $E([\al,\xi])$.  The image $dW$ of $W$ under the
differential $d: \Mm^5\to \Mm^{6}$ is dual to the above Whitehead product
map $[\cdot,\cdot]$.  Hence $dW = A(X-Y)$.  

Moreover, if
$\Aa$ denotes the vector space of indecomposable
(or primitive) 
elements in $\Mm$ and  $\Mm_{\ge k}$ is the subalgebra generated by
  $\Aa^{\otimes k}$ 
the differential $d: \Mm \to \Mm/\Mm_{\ge k+2}$ is determined 
by the nonvanishing Whitehead products of order $\le k$.  In the case of
$BG_\la$ for $0< \la \le 1$ it is only the first order, i.e. the usual, Whitehead
products, that can be nonzero, since these account for the whole of
$\pi_*(G)$.  It follows that $H^*(BG)$ has only the one relation
$A(X-Y) = 0$.\QED

\begin{rmk}\rm   Another way of seeing the relation between
 Samelson (or Whitehead) products  and $H^*(BG)$ is to look at the 
cohomology spectral sequence of the fibration
$G\to EG\to BG$.  In general, the relation is complicated, and the more
natural spectral sequence to look at is one constructed by Quillen: see
Allday~\cite{Al}.  However, in the case considered above it is not too hard
to work out what happens since only first order products are involved.
The differential
$$ 
d_2^{0,q}: H^0(BG)\otimes H^q(G)\to
H^{2}(BG)\otimes H^{q-1}(G) 
$$
 takes $a$ to $A$, and $x,y$ to $0$, and, when $q=4$, its restriction to the
 the primitive part of $H^4(G)$ is dual to the 
Samelson product
$$
H_2(BG)\otimes H_3(G) \cong \pi_1(G)\otimes \pi_3(G)\otimes \Q \to \pi_4(G)\otimes \Q.
$$
 Thus 
$$
d_2^{0,4}(w) = A\otimes (x-y) = d_2(a\cup(x-y)),
$$
so that $$
E_3^{0,4} \cong E_4^{0,4} \cong \Q \widehat{w}.
$$
  Similarly,
$$
d_4^{0,4} : H^0(BG)\otimes H^4(G)\to
H^{4}(BG)\otimes H^{1}(G)
$$
is dual to the Samelson product and so takes the generator
$\widehat{w}$ of $E_4^{0,4}$ to the  element $a(X-Y)$ of
$$
E_4^{4,1} = \Ker(d_2^{4,1}) \subset H^4(BG)\otimes H^1(G).
$$
Thus $a(X-Y)$ must lie in the kernel of $d_2^{4,1}$, i.e. $A(X-Y) = 0$ in $H^*(BG)$.
\end{rmk}

Next, let us consider the case $1 < \la \le 2$, i.e. $\ell = 2$.  
We will discuss this case in detail since it is a paradigm for the others.
Now, the first order Samelson products $[\al,\xi], [\al,\eta]$ both
vanish since  $[\al, \xi+\eta] = 0$ as before, and the classes $\al$ and $ \xi +
4\eta $ have  commuting representatives in
$K_2$.  Hence all first order Samelson products  vanish.  
 Therefore, the second order
Samelson product 
$$ 
[\al,\xi,\xi] \in \pi_8(G)\otimes \Q
$$ 
is defined.  As explained in~\cite{Al}, in general $[\al,\xi,\xi] $ should be
considered as a coset in $\pi_8(G)\otimes \Q$, but this is reduced to a single
element because the relevant homotopy groups vanish.  More precisely,
given two maps $\be_i: S^{p_i} \to X$, $i= 1,2$, whose Samelson product
$$
[\be_1,\be_2]: S^{p_1+p_2}\to X
$$
is null-homotopic, let us write $C[\be_1,\be_2]$ for the chain given by a
particular choice of null-homotopy $D^{p_1 + p_2 + 1}\to X$. (This is unique
up to homotopy.) Then, Allday shows in~\cite{A}~\S2 that the image of
$[\al,\xi,\xi] $ in $H_8(G)$ is represented by 
$$
[ \al,\xi,\xi] = [C[\al,\xi], \xi] + [C[\xi,\xi], \al] + [C[\xi,\al], \xi]
 = 2[C[\al,\xi], \xi] + [C[\xi,\xi], \al],
$$
since $[\al,\xi] = [\xi, \al]$.  This is a cycle because of the Jacobi identity. 
Its suspension is the obstruction to extending the map
$$
E(\al)\vee E(\xi)\vee E(\xi)\;:\quad S^2\vee S^4\vee S^4 
\longrightarrow BG
$$
to the product $S^2\times S^4\times S^4$.

\begin{lemma}\label{le:sa2} When $\ell = 2$, 
\begin{itemize}
\item[(i)]
$[ \al,\xi+\eta,\xi+\eta] = [ \al,\xi+4\eta,\xi+4\eta]=0$. 
\item[(ii)]  $w_2([\al,\xi,\xi]) \ne 0$. 
\end{itemize} \end{lemma}
\proof{} First observe that when $k = 1$ or $2$, 
one can 
choose commuting representatives for $\al$ and $\xi + k^2\eta$ in $K_k$. 
Moreover, one can also choose a representative for  $C[\xi+k^2\eta,\xi+k^2\eta]$
that lies in the $\SO(3)$ factor of $K_k$ and so commutes with $\al$. 
This proves (i).  

Now consider (ii).
 As in Lemma~\ref{le:sa1}, let us look at
the image of this cycle in $U_0$. Since $[\xi,\xi] = 0$ in $H_*(K_0)$ we can
choose $C[\xi,\xi]$ to lie in $K_0$. Therefore, in $U_0 = G/K_0$ the class $
[ \al,\xi,\xi]$ is represented by
$$
2\xi * C[\al,\xi] + C[\xi,\xi] * \al.
$$
Note that $G$ acts on $\Jj = \Jj_\la$ by $J\mapsto \phi_*(J)$, and we
 may choose the neighborhood $\Nn(U_1)$ to be invariant under 
the induced action of the compact group $K_0$.  Further, given a $K_0$ action
$K_0\times V\to V$, we write $*$ to denote the induced product 
$$
*:\quad H_i(K_0)\otimes
H_j(V)\to H_{i+j}(V).
$$
 If we  choose $J\in \Nn(U_1)
- U_1$ as in Lemma~\ref{le:sa1}, the chains $ C[\xi,\xi]$,  $ \al,$ and $  C[\xi,\xi] *
\al$ all have representatives in $\Nn(U_1) - U_1$.   In particular, $ C[\xi,\xi] * \al$ 
can be considered as an element of $H_8(\Nn(U_1), \p \Nn(U_1))$ with boundary
$-2\xi*[\al,\xi]\in H_7(\p\Nn(U_1))$.   But 
$$
H_8(\Nn(U_1), \p \Nn(U_1))\cong H_4(U_1) = 0.
$$
Thus we can replace the chain $ C[\xi,\xi] * \al$ (which does not intersect
$U_1$) by any  chain with the same boundary.
But there is an obvious chain in $\Nn(U_1)$ with boundary $[\al,\xi] =
\xi*\al$, namely, $\xi*D$ where $D$ is a $2$-disc in $\Nn(U_1)$ that is
transverse to the stratum $U_1$.  Thus 
$ [ \al,\xi,\xi]$ is represented by
$$
\xi *(C[\al,\xi] - 2\xi*D),
$$
where $C[\al,\xi] - 2\xi*D$ is a $5$-cycle in $U_{01}$.
Observe that $D$ has non trivial intersection with $U_1$ while 
$C[\al,\xi]$ lies in $U_0$.  Hence the cycle 
$C[\al,\xi] - 2\xi*D$ has
nontrivial linking number with $U_2$.  In other words, its image under
the map $\tau: H^5(U_{01})\to H^0(U_2)$ of $(MV_2)$ is nonzero.  The
desired conclusion now follows from Definition~\ref{def:w}.\QED

\begin{cor}\label{cor:BG2}When $1< \la\le 2$, $H^*(BG_\la) =
\Q[A,X,Y]/\{A(X-Y)(4X-Y) = 0\}$, where $A$ has dimension $2$ and $X,Y$ 
have dimension $ 4$. \end{cor}
\proof{}  Again, the relations in $H^*(BG)$ are dual to the nonzero
higher Whitehead products in $\pi_*(BG)$: see~\cite{AA}~Thm~5.4.  
Since $[\al,\al] = 0$, the three elements $\al, \xi, \eta$ 
are transgressive and give rise to nonzero elements $A,
X,Y$ of $H^*(BG_\la)$. 
Moreover, since only one of the groups 
$\pi_*(BG)\otimes \Q = \pi_{*-1}G\otimes \Q$ for $* > 4$ is
nonzero, there can be at most one  relation between these generators.  In fact, the
question here is to decide whether the new element $w_\ell$ is transgressive
in the spectral sequence of the fibration $G \to EG \to BG$, i.e.
gives rise to a new generator of $H^*(BG)$, or whether it gives rise to a relation
between the existing generators.   The previous lemma tells us that the latter holds.
Thus there is exactly one nontrivial relation $F = 0$ in $H^*(BG)$.  Since it
corresponds to a  second order Whitehead product, $F$ is homogeneous of order $3$
in the variables $A,X,Y$, and we can think  of it as a symmetric trilinear function on
the vector space  $V$ that is spanned over $\Q$ by a basis $\{e_\al, e_\xi,
e_\eta\}$ dual
 to $\{A,X,Y\}$.  
The
vanishing results in Lemma~\ref{le:htpyK0} (i) tell us that 
$$
F(e_\al, e_\xi+e_\eta, e_\xi+e_\eta) = 0,\quad F(e_\al, e_\xi+4e_\eta,
e_\xi+4e_\eta ) = 0.
$$
Since $F(e_\al,e_\xi, e_\xi)\ne 0$, $F$ has to be a nonzero multiple of $A(X - Y) (4X
- Y)$, as claimed.\QED

When $\ell > 2$ the argument is similar.  We first  explain
higher order Whitehead products.   Suppose given maps
$\be_i:S^{p_i} \to X$ such that 
$$
 \be_1\vee\dots\vee \be_k\;:\quad S^{p_1}\vee \dots \vee S^{p_k}
\longrightarrow X
$$
extends to a map $f$ of  the $(p-1)$-skeleton $T_I$ of $P_I = S^{p_1}\times
\dots\times S^{p_k}$ into $X$, where $p = \sum_i p_i$ and $I = (p_1,\dots, p_k)$. 
Then the usual higher Whitehead product $[\be_1,\dots\be_k]$ is defined to be the
set of all obstructions $Ob(f)$ to extending $f$ over $P_I$, where $f$ ranges over all
possible extensions to $T_I$.  Thus this product is a coset in
$\pi_{p-1}(X)\otimes \Q$.  In our situation the space $X$ corresponds to $BG$ and  we must
consider the analogous higher order products in $G$.  These  are  also
called   higher Whitehead products (though it might be more logical to call them
higher Samelson products.)  As usual the $p$th order products are defined only
when all relevant lower order ones vanish. 
According to~\cite{Al2} \S3 they can be defined inductively as follows.
Again, we will be somewhat sloppy with notation, and will consider the formulas
below to define representatives of the given classes, as well as  the classes
themselves. 

Suppose given representatives $\xi_i$ of elements $\xi_i\in \pi_{n_i}(G)\otimes
\Q$, and let $I$ be a proper subset $\{j_0,\dots j_k\}$ of $I_n =  \{0,1,2,\dots, n\}$. 
Then the $k$th order Whitehead product 
$$
\nu_I = [ \xi_{j_0},\dots, \xi_{j_k} ] \in \pi_{N_I}(G)\otimes \Q
$$
 of the elements
$\xi_{j_0},\dots, \xi_{j_k}$ has degree $
N_I = k + \sum_{j\in I} n_j $ if it is defined.  
If $\nu_I$ is null-homotopic, then the corresponding map from the $N_I$-skeleton
$T_I$ of  the product of spheres 
$$
P_I = S^{n_1 + 1}\times \dots\times S^{n_k+1}
$$
 extends to $P_I$ and, as above, we
write $C\nu_I$ for this extension.  We should think of $C\nu_I$ as a singular
chain with boundary $\nu_I$.  (Usually it
depends on choices, though in our situation it will not.) We extend this definition
to singleton sets $I=\{j\}$, by setting $C\nu_I = \xi_j$ in this case. 

Now suppose that these chains $C\nu_I$ are defined for all proper
subsets $I$ of $I_n = \{0,1,\dots, n\}$.  In this situation Allday shows that a
representative for the $n$th order product is given by a  formula of the type $$
\nu_{I_n}  = [ \xi_{0},\dots, \xi_{n} ] = \sum_{I\subset I_n} \pm [C\nu_I,
C\nu_{I_n - I}] $$ where the sign depends on the shuffle $(I, I_n - I)$.\footnote 
{ In fact, Allday interprets this formula  in the minimal model of $X=BG$, showing
that it defines a
 primitive element there.  Therefore it
is spherical and   can be considered as an element in $\pi_*(BG)\otimes \Q =
\pi_{*-1}(G)\otimes \Q$.}

In our situation, all intermediate homotopy groups vanish and so
  the $n$th order
product $\zeta_n= [ \al, \xi,\dots, \xi ]$ represents a single element of
$\pi_{4n}(G)\otimes \Q$.  Moreover, it is determined in terms of the $\zeta_k$ for $0 \le k < n$
and in terms of the $p$th order products 
$\xi_p = [  \xi,\dots, \xi ]$ which we can consider as (null homologous)
cycles in $K_0$.  (Here we take $\zeta_0 = \al$.)  Moreover,
$$
\zeta_n =  [ \al, \xi,\dots, \xi ]
=\ell [\xi, C\zeta_{n-1}]  + \mbox{ lower order terms}
$$
where by lower order terms we mean terms of the form $[C\xi_p, C\zeta_{n-p}]$
with $p > 1$. 

In the next lemma we think of $\xi_p, C\xi_p$ as chains in $K_0$
that, as in Lemma~\ref{le:sa2}, act via the Pontrjagin product $*$ on chains in
$\Jj_\la$.  In particular, we denote by $S$ the link of $U_\ell$ in
$\Jj_\la$. Thus $S\cong S^{4\ell -3}$ and $\xi*S$ is a cycle in $U_{0\dots \ell-1}$.
By Remark~\ref{rmk:w}, this cycle generates $H_{4\ell}(U_{0\dots \ell-1})\cong
H_{4\ell}(U_0)$.  Hence the next lemma implies immediately that
$w_\ell(\zeta_\ell)\ne 0$.

 \begin{lemma}  If $\ell-1 < \la \le \ell$, the $\ell$th order
Whitehead product $\zeta_\ell\in \pi_{4\ell}(G_\la)\otimes\Q$ is defined and 
its image in $U_{0\dots \ell-1}\subset \Jj_\la$ may be represented by
a nonzero multiple of the cycle $\xi*S$ described above.
\end{lemma}
\proof{}  We prove this by induction on $\ell$.  It is true when $\ell = 1,2$ by
our previous work.  To prove it for $\ell > 2$ we start, as in Lemma~\ref{le:sa2},
 from the equation
$$
 \zeta_\ell =  [ \al, \xi,\dots, \xi ]
= \ell\, [\xi, C\zeta_{\ell-1}]  + \mbox{ lower order terms}.
$$
Since  $\xi_p$ and $C\xi_p$ lie in $K_0$ the image of 
$\zeta_\ell$ in $U_0 =
G/K_0$ has the form
\begin{equation}\label{eq:ze}
\ell \xi*C\zeta_{\ell-1} + \mbox{terms of the form } C\xi_p * C\zeta_{\ell-p}.
\end{equation}
Again, the boundary of $\xi*C\zeta_{\ell-1}$ is $\xi*\zeta_{\ell-1}$.
By the inductive hypothesis, the
cycle $\zeta_{\ell-1}$ can be represented in a
 deleted neighborhood $\Nn(U_{\ell -1}) - U_{\ell -1}$ of
$U_{\ell -1}$.  Moreover, if we choose this neighborhood to be 
$K_0$-invariant, we can  represent  $\xi*\zeta_{\ell-1}$ 
and all the cycles $C\xi_p * C\zeta_{\ell-p}, p> 1,$ in the above formula there as
well. Note, however, that $C\zeta_{\ell-1}$ cannot be represented there
since, by the inductive hypothesis $\zeta_{\ell-1} = \xi*S$, which  is
clearly  nonzero
in $H_*(\p\Nn(U_{\ell -1}))$.  On the other hand, $S$ is bounded by
a ball $D$ in the full neighborhood
$\Nn(U_{\ell-1})$ so that $\zeta_{\ell-1}$ is the boundary of $\xi*D$ in 
$\Nn(U_{\ell-1})$.  Therefore, as in Lemma~\ref{le:sa2},
we can replace the sum of terms $C\xi_p * C\zeta_{\ell-p}$ in equation~(\ref{eq:ze})
by  a suitable multiple of  $\xi*D$.  The result now follows as before.\QED

\begin{cor} When $\ell -1< \la\le \ell$, 
$$
H^*(BG_\la) =
\Q[A,X,Y]/\{A(X-Y)(4X-Y)\dots(\ell^2X - Y) = 0\},
$$
 where $A$ has dimension $2$ and
$X,Y$  have dimension $ 4$. 
\end{cor}
\proof{}  By Lemma~\ref{le:htpyK0}
all $\ell$th order Whitehead products of the form
$$
[\al, \xi+ k^2\eta, \dots, \xi+ k^2\eta], 1 \le k \le \ell,
$$
vanish in $\pi_*(G_\la)$.  The rest of the proof is much the same as
the proof of Corollary~\ref{cor:BG2} and is left to the reader.\QED

This proves Theorem~\ref{thm:bg}.   The proof of the corresponding
result for $\cp\#\,\ocp$ will, as usual, be left to the reader.

\end{document}